\documentclass[]{iopart}
\usepackage{iopams,graphicx,hyperref,hhline}
\usepackage[normalem]{ulem}
\usepackage{subcaption}
\usepackage[usenames,dvipsnames]{color}
\usepackage{mathrsfs}
\DeclareSymbolFontAlphabet{\mathrsfs}{rsfs}
 
\newcommand{\tfrac}[2]{\textstyle \frac{#1}{#2}} 
\newcommand{\R}{\mathbb{R}}
\renewcommand{\Re}{\mathrm{Re}}
\renewcommand{\Im}{\mathrm{Im}}
\newcommand{\half}{\tfrac{1}{2}}

\newcommand{\Lie}{\mathcal{L}}
\newcommand{\const}{\mathrm{const}}
\newcommand{\Scri}{$\mathrsfs{I}^+\,$} 
\renewcommand{\tr}{\tilde r} 
\newcommand{\htr}{h_{\tilde r}} 
\newcommand{\Ntr}{N_{\tilde r}}
\newtheorem{conj}{Conjecture}

\begin{document}

\title[Hyperboloidal numerical method for multidimensional nonlinear wave equations]{A hyperboloidal method for numerical simulations of multidimensional nonlinear wave equations: nonlinear tails}
\author{Oliver Rinne}
\address{HTW Berlin -- University of Applied Sciences, Faculty 4, Treskowallee 8, 10318 Berlin, Germany}
\ead{oliver.rinne@htw-berlin.de}

\begin{abstract}
  We consider the scalar wave equation with power nonlinearity in $n+1$ dimensions. 
  Unlike most previous numerical studies, we go beyond the radial case and do not assume any symmetries for $n=3$, and we only impose an SO$(n-1)$ symmetry in higher dimensions.
  Our method is based on a hyperboloidal foliation of Minkowski spacetime and conformal compactification.
  We focus on the late-time power-law decay (tails) of the solutions and compute decay exponents for different spherical harmonic modes, for subcritical, critical and supercritical, focusing and defocusing nonlinear wave equations.
\end{abstract}


\section{Introduction}

This paper is concerned with the nonlinear wave equation (NLW)
\begin{equation} 
  \label{e:nlw1}
  \Box \Phi := -\partial_t^2 \Phi + \Delta \Phi =  \mu |\Phi|^{p-1} \Phi, \quad \Phi : \R \times \R^n \to \R,
\end{equation}
where $p>1$ and $\mu=\pm 1$, with $\mu=-1$ referred to as the focusing NLW and $\mu=1$ as the defocusing NLW. 

This equation serves as a model for various nonlinear wavelike equations arising e.g. in fluid dynamics, optics, acoustics, plasma physics, general relativity and quantum field theory.
Related equations include the nonlinear Schr\"odinger equation, the Korteweg-de Vries equation, the Klein-Gordon equation, Yang-Mills equations and wave map equations.

The NLW \eref{e:nlw1} is invariant under the rescaling
\begin{equation} \label{e:rescaling}
  \Phi(t,x) \to \lambda^\frac{2}{p-1} \Phi(\lambda t, \lambda x).
\end{equation}
The equation possesses a conserved energy 
\begin{equation} \label{e:stdenergy}
  E(\Phi, \partial_t \Phi) = \int_{\R^n} \left( \half (\partial_t \Phi)^2 + \half \Vert \nabla \Phi \Vert^2 + \frac{\mu}{p+1} |\Phi|^{p+1} \right) \rmd x.
\end{equation}
The energy is invariant under the rescaling \eref{e:rescaling} iff
\begin{equation} \label{e:pcrit}
  p = \frac{n+2}{n-2} =: p_\mathrm{crit},
\end{equation}
in which case the NLW is called energy-critical.
For $p<p_\mathrm{crit}$ the NLW is subcritical, for $p>p_\mathrm{crit}$ supercritical.

Solutions to the NLW can show rich behaviour due to the interplay of the dispersive wave operator and the nonlinear term. 
While small initial data will generally scatter, i.e. approach a solution to the \emph{linear} wave equation, in the focusing case large initial data will cause blow-up of the solution. 
The threshold behaviour between scattering and blow-up is particularly interesting and may involve a universal attractor \cite{Bizon2004,Bizon2009}.
Futhermore, in the critical case stable finite-energy solitons exist that may prevent solutions from scattering \cite{Aubin1976,Talenti1976,Kenig2008}.

In this paper we will mainly explore the late-time behaviour of solutions that scatter.
In \cite{Szpak2009} it was proved using perturbative methods that in $n=3$ spatial dimensions, spherically symmetric solutions to \eref{e:nlw1} decay as $t^{-p+1}$ for $p\geqslant 3$.
This is often referred to as a power-law \emph{tail}.
Here we will study such tails numerically beyond the spherically symmetric case.

\begin{figure} 
  \begin{subfigure}[b]{0.5\textwidth}
    \includegraphics[width=\textwidth]{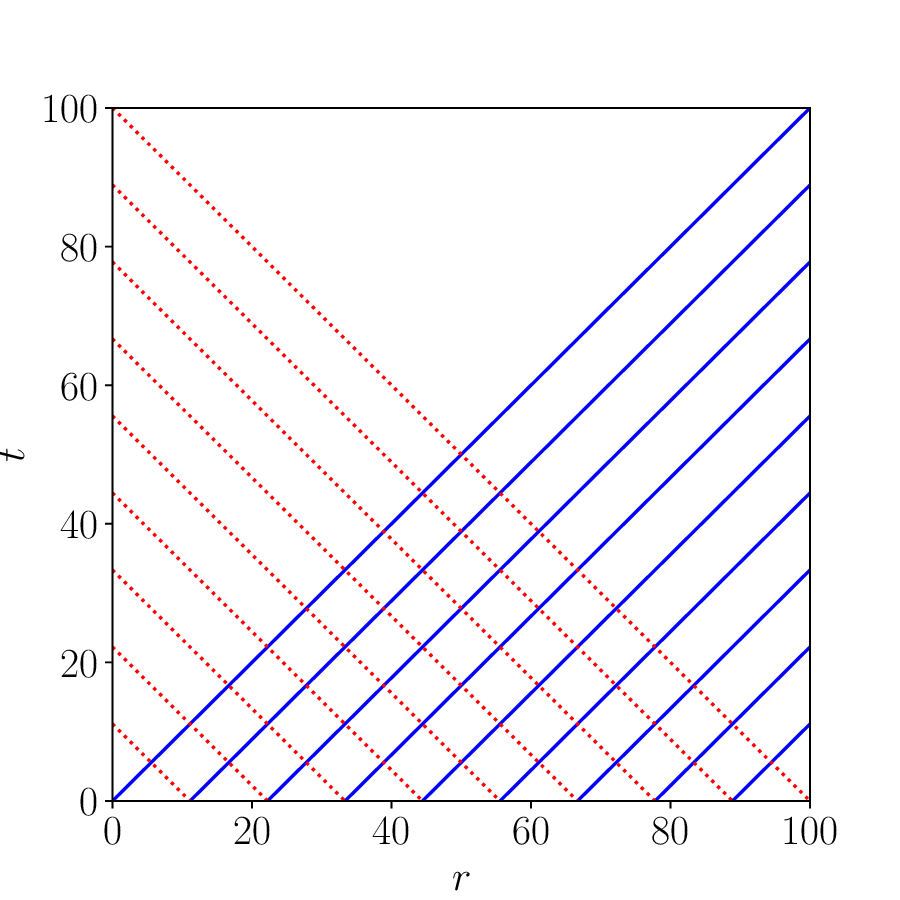}
    \caption{\label{f:char_std}}
  \end{subfigure}
  \begin{subfigure}[b]{0.5\textwidth}
    \includegraphics[width=\textwidth]{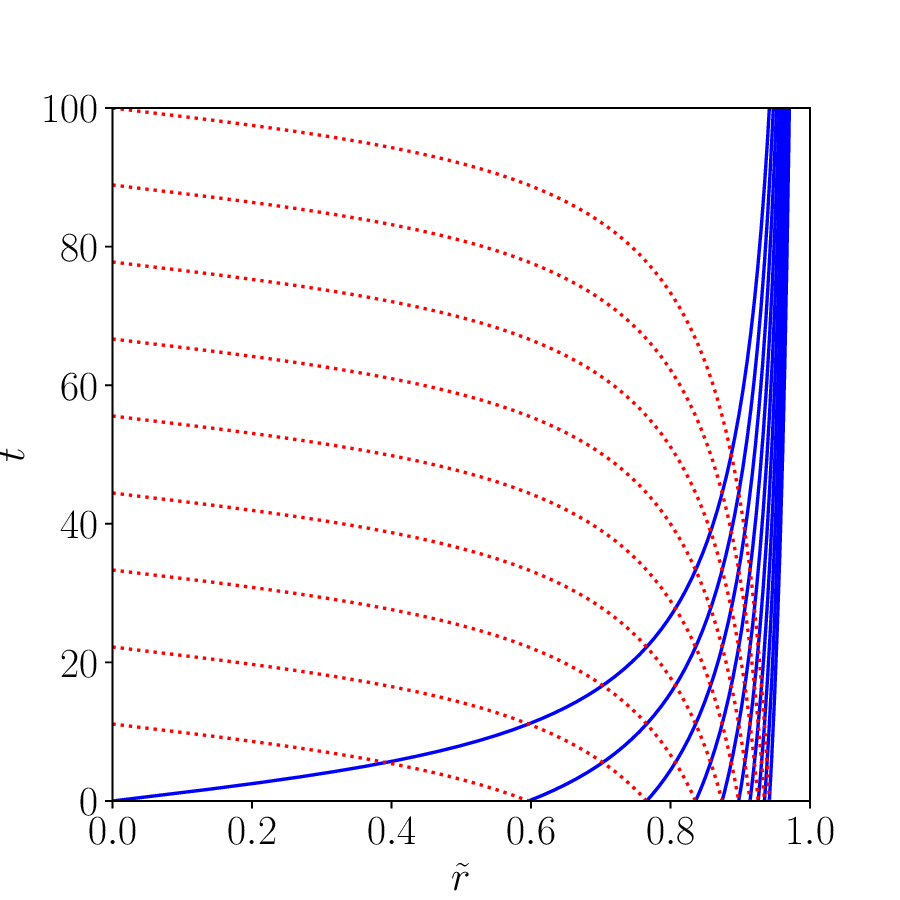}
    \caption{\label{f:char_spt_cpt}}
  \end{subfigure}\\
  \begin{subfigure}[b]{0.5\textwidth}
    \includegraphics[width=\textwidth]{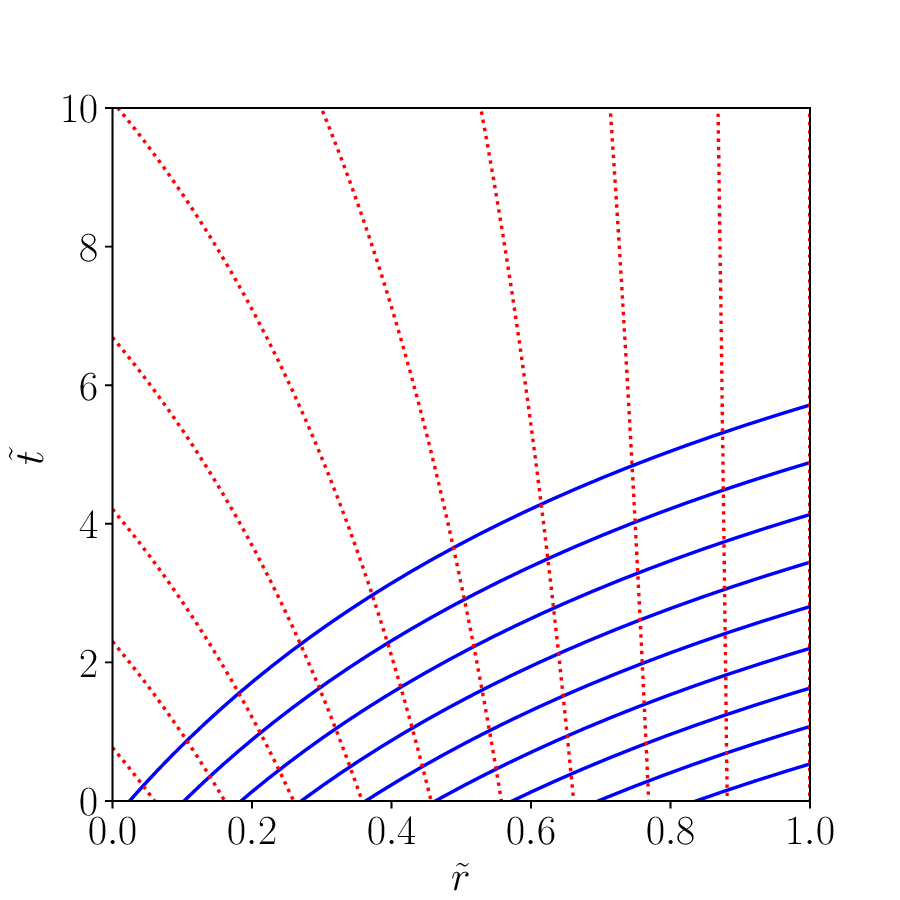}
    \caption{\label{f:char_hyp_cpt}}
  \end{subfigure}
  \begin{subfigure}[b]{0.5\textwidth}
    \includegraphics[width=\textwidth]{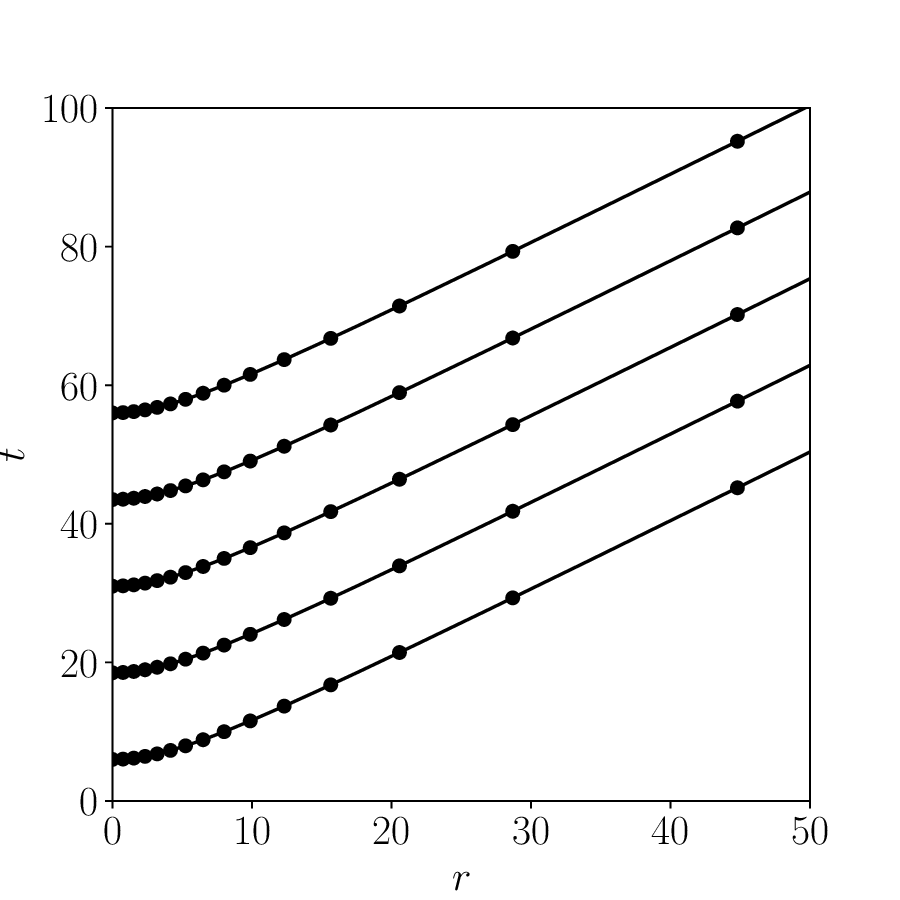}
    \caption{\label{f:hyp_slices}}
  \end{subfigure}
  \caption{\label{f:hyperboloidal} 
    Illustration of the hyperboloidal method. 
    (a) Outgoing (solid blue) and ingoing (dotted red) characteristics $t=\pm r + \const$ of the wave equation.
    (b) The same characteristics plotted against the compactified coordinate $\tilde r$ defined in \eref{e:rtilde_intro}.
    (c) Characteristics in terms of $\tilde r$ and the hyperboloidal time coordinate $\tilde t$ defined in \eref{e:ttilde_intro}.
    (d) A few hyperboloidal slices $\tilde t = \const$ (evenly spaced in $\tilde t$), plotted in $(r,t)$ coordinates. The dots indicate a grid that is uniform in the compactified coordinate $\tilde r$.
    We have chosen $a=6$ in \eref{e:rtilde_intro} and \eref{e:ttilde_intro}, which corresponds to dimension $n=3$ and mean curvature constant $C=0.5$ in section \ref{s:metric}.
  }
\end{figure}

In the standard approach to solving wavelike equations numerically, applied e.g. in \cite{Strauss1978,Colliander2010,Donninger2011,Murphy2020}, the spatial domain is taken 
to be a large ball.
Boundary conditions must be imposed at the surface of this ball in order to obtain a well-posed initial-boundary value problem.
Typically a homogeneous Dirichlet boundary condition is used. 
This causes spurious reflections once the outgoing waves reach the boundary so that the numerical solution can only be trusted up to a certain time.
Mapping the radial coordinate $r$ from $(0,\infty)$ to a finite interval, e.g.
  \begin{equation} \label{e:rtilde_intro}
    r = \frac{2a\tilde r}{1-\tilde r^2}, \quad \tilde r \in (0,1), \quad a=\const,
  \end{equation}
and discretising the compactified coordinate $\tilde r$ is not a reliable approach either because the wavelength w.r.t. $\tilde r$ decreases towards zero as the waves travel outwards, and they ultimately fail to be resolved numerically.
This is illustrated in figure \ref{f:hyperboloidal}, where we plot the characteristics of the wave equation in terms of $r$ (figure \ref{f:char_std}) and $\tilde r$ (figure \ref{f:char_spt_cpt}).

A different approach is to introduce a new time coordinate $\tilde t$ such that the slices $\tilde t = \const$ become asymptotically characteristic, e.g.
\begin{eqnarray} \label{e:ttilde_intro}
  \tilde t = t - \sqrt{a^2 + r^2}
\end{eqnarray}
with a constant $a$. 
Such \emph{hyperboloidal} slices approach \emph{future null infinity} \Scri as $r\to\infty$, the asymptotic region where outgoing null geodesics end.
By means of a suitable compactification, \Scri can be brought to a finite radius on the numerical grid.
No boundary conditions are needed there because all characteristics point towards the exterior of the domain, and the infinite ``blueshift'' of outgoing waves is eliminated (figure \ref{f:char_hyp_cpt}).
In figure \ref{f:hyp_slices} we plot a few hyperboloidal slices in the $(r,t)$ plane, along with a numerical grid that is chosen to be uniform in $\tilde r$.

The hyperboloidal method originated in general relativity, see \cite{Frauendiener2004} for a comprehensive review article.
It has been applied to nonlinear wave equations in spherical symmetry \cite{Bizon2009}, spherically symmetric linear scalar wave and Yang-Mills equations in fixed Schwarzschild spacetime \cite{Zenginoglu2008} and coupled to the Einstein equations \cite{Rinne2013}, and linear scalar wave equations without symmetries in Kerr spacetime \cite{Racz2011}.
The cubic ($p=3$) focusing ($\mu=-1$) NLW \eref{e:nlw1} in $n=3$ spatial dimensions without symmetries was evolved in \cite{Zenginoglu2010,Zenginoglu2011}.

As far as we are aware, the present paper is the first extensive numerical study of the NLW in higher dimensions beyond spherical symmetry. 
We do not assume any symmetries in $n=3$ spatial dimensions. 
In higher dimensions, we impose an SO$(n-1)$ symmetry so that there is one effective angular coordinate. 

Our numerical method is similar to that in \cite{Racz2011,Csizmadia2012} in that it combines a finite-difference method in the radial direction with a spectral method in the angular directions.
However, due to the nonlinear terms in our wave equation, we employ a pseudo-spectral collocation method.

This paper is organised as follows.
In section \ref{s:formulation} we derive the form of the wave equation in the framework of the hyperboloidal method.
In section \ref{s:nummeth} we describe the numerical methods we use to solve this equation.
We present our numerical results in section \ref{s:results}: a convergence test against exact linear solutions, a test of the energy balance on hyperboloidal slices, and power-law tails.
The results are summarised and discussed in section \ref{s:conclusion}.
Exact solutions to the linear wave equation are derived in \ref{s:exact}.


\section{Formulation} 
\label{s:formulation}

In this section we derive the form of the wave equation we will be solving numerically.
In section \ref{s:metric} we foliate Minkowski spacetime into hyperboloidal slices of constant mean curvature and apply a conformal transformation to the metric. 
We work out the wave equation in terms of the conformally rescaled quantities in section \ref{s:waveeq} and write it in first-order in time, second-order in space form.
Finally we derive the energy balance that the scalar field obeys on our hyperboloidal foliation in section \ref{s:energy_theory}.


\subsection{Hyperboloidal foliation of Minkowski spacetime}
\label{s:metric}
 
We start with the $(n+1)$-dimensional Minkowski metric in spherical polar coordinates,
\begin{equation}
  \eta = -\rmd t^2 + \rmd r^2 + r^2 \sigma^{(n-1)},
  \quad t\in\R, \quad r\in(0,\infty).
\end{equation}
Here $\sigma^{(n-1)}$ is the standard round metric on $S^{n-1}$, the $(n-1)$-dimensional unit sphere. 
Explicitly, we have 
\begin{eqnarray}
  \sigma^{(1)} = \rmd \varphi^2, \\
  \sigma^{(2)} = \rmd \theta^2 + \sin^2\theta \, \rmd\varphi^2 , 
\end{eqnarray}   
where $\theta \in [0,\pi]$ and $\varphi \in [0,2\pi)$.
More generally for $n>3$,
\begin{equation} 
  \sigma^{(n-1)} = \rmd \theta_1^2 + \sin^2\theta_1 \left( \rmd \theta_2^2 + \sin^2\theta_2 (\ldots (\rmd \theta_{n-2}^2 + \sin^2 \theta_{n-2} \rmd\varphi^2 ) )\right)
\end{equation}
with $\theta_1, \ldots, \theta_{n-2} \in [0,\pi]$ and $\varphi \in [0,2\pi)$.

We introduce a new time coordinate 
\begin{equation}
  \label{e:hyperboloid}
  \tilde t = t - \sqrt{a^2 + r^2}
\end{equation}
with a constant $a = n/C > 0$. 
The hypersurfaces $\tilde t=\const$ are hyperboloids and $C$ is their constant mean curvature (up to a sign as discussed further below).
In the new coordinates the Minkowski metric reads
\begin{equation}
  \eta = -\rmd \tilde t^2 - \frac{2r}{\sqrt{a^2 + r^2}} \, \rmd \tilde t \, \rmd r + \frac{a^2}{a^2 + r^2} \, \rmd r^2 + r^2 \sigma^{(n-1)}. 
\end{equation}
Let us also introduce a new radial coordinate $\tr$ depending on $r$ only.
We can arrange for the spatial metric to be conformally flat if we demand 
\begin{equation} \label{rdiff}
  \frac{a^2}{a^2 + r^2} \rmd r^2 = \frac{r^2}{\tr^2} \rmd \tr^2 \; \Leftrightarrow \; \frac{\rmd r}{\rmd \tr} = \frac{r \sqrt{a^2 + r^2}}{\tr a}.
\end{equation}
This ordinary differential equation can be solved explicitly.
If we impose the boundary condition $\tr \to 1$ as $r\to \infty$, we obtain 
\begin{equation}
  \label{e:rofr_}
  r = \frac{2a\tr}{1-\tr^2} = \frac{2n\tr}{C(1-\tr^2)}.
\end{equation}
From \eref{e:hyperboloid} we observe that on a slice $\tilde t=\const$, as the conformal boundary $\tr\to 1$ is approached, we have  
\begin{equation}
  t, r \to \infty, \quad t - r \to \tilde t = \const,
\end{equation}
so the slice becomes asymptotically characteristic and $\tr=1$ corresponds to \Scri.

The metric now takes the form 
\begin{eqnarray}
  \eta &=& -\rmd \tilde t^2 - \frac{2 r^2}{\tr a} \rmd \tilde t \, \rmd \tr + \frac{r^2}{\tr^2} \left( \rmd \tr^2 + \tr^2 \sigma^{(n-1)} \right)\nonumber\\
  &=& \Omega^{-2} \left[ - \tilde \alpha^2 \rmd \tilde t^2 + (\rmd \tr + \tilde \beta^{\tr} \rmd \tilde t)^2 + \tr^2 \sigma^{(n-1)} \right] =: \Omega^{-2} \tilde \eta. 
\end{eqnarray}
Here we have introduced the conformal factor 
\begin{equation} \label{e:Omega}
  \Omega = \frac{\tr}{r} = \frac{C}{2n} (1-\tr^2),
\end{equation}
the conformal shift vector $\tilde \beta$, which only has a radial component 
\begin{equation} \label{e:shift}
  \tilde \beta^{\tr} = -\frac{\tr}{a} = -\frac{C}{n}\tr,
\end{equation}
and the conformal lapse function 
\begin{equation} \label{e:lapse}
  \tilde \alpha = \frac{C}{2n} (1+\tr^2)
\end{equation}
so that $\tilde \alpha^2 - (\tilde \beta^{\tr})^2 = \Omega^2$. 

From now on we will work entirely with the conformal metric $\tilde \eta$.
Its induced metric $\tilde \gamma$ on the $\tilde t=\const$ slices is flat,
\begin{equation}
  \tilde \gamma = \rmd \tr^2 + \tr^2 \sigma^{(n-1)}.
\end{equation}
The conformal extrinsic curvature of the $\tilde t=\const$ slices is\footnote{
  We use the sign convention $\Lie_{\tilde \nu} \tilde \gamma_{ab} = -2 \tilde K_{ab}$.
  For the $\tilde t=\const$ slices to intersect \Scri we need the constant $C>0$ so that the shift \eref{e:shift} is negative at $\tr=1$.
  Spacetime indices $a,b$ are raised and lowered using the conformal metric $\tilde \eta$, and repeated indices are summed over.
} 
\begin{equation} \label{e:confKab}
  \tilde K_{ab} = \frac{\tilde K}{n} \tilde \gamma_{ab}
\end{equation}
with the conformal mean curvature given by
\begin{equation}
  \tilde K = -\tilde \alpha^{-1} C = -\frac{2n}{\tr^2+1}. 
\end{equation}
Its Lie derivative along the conformal future-directed timelike unit normal $\tilde \nu$ to the $\tilde t=\const$ slices is 
\begin{equation}
  \Lie_{\tilde \nu} \tilde K = -\tilde \alpha^{-1} \tilde \beta^{\tr} \partial_{\tr} \tilde K = \frac{8n \tr^2}{(\tr^2+1)^3}.
\end{equation}

Combining the contracted Gauss and Ricci equations, we can relate the Ricci scalar $\tilde R$ of $\tilde \eta$ to the Ricci scalar $\tilde R_{\tilde \gamma}$ of $\tilde \gamma$:
\begin{equation} \label{e:gaussricci}
  \tilde R = \tilde R_{\tilde \gamma}- 2 \tilde \alpha^{-1} \tilde \Delta \tilde \alpha - 2 \Lie_{\tilde \nu} \tilde K + \tilde K^2 + \tilde K_{ab} \tilde K^{ab}.
\end{equation}
Here 
\begin{equation}
  \tilde \Delta = \partial_{\tr}^2 + (n-1) {\tr}^{-1} \partial_{\tr} + \tr^{-2} \mathring{\Delta}^{(n-1)}
\end{equation}
is the standard flat Laplacian of the flat conformal metric $\gamma$ and
$\mathring{\Delta}^{(n-1)}$ is the Laplace-Beltrami operator on the sphere $S^{n-1}$.
Explicitly, we have
\begin{eqnarray}
  \mathring{\Delta}^{(1)} = \partial_\varphi^2, \\
  \mathring{\Delta}^{(2)} = \partial_\theta^2 + \cot\theta \, \partial_\theta + \frac{1}{\sin^2\theta}\partial_\varphi^2.
\end{eqnarray}

In our numerical implementation for higher dimensions $n>3$, we shall impose an SO$(n-1)$ symmetry so that all but one of the angular coordinates are suppressed.
The angular Laplacian then reduces to
\begin{equation}
  \mathring{\Delta}^{(n-1)} = \partial_\theta^2 + (n-2)\cot\theta \, \partial_\theta.
\end{equation}

Using the fact that $\tilde R_{\tilde \gamma} = 0$ and the form of the conformal extrinsic curvature \eref{e:confKab}, we obtain from \eref{e:gaussricci}
\begin{equation}
  \tilde R = -2 \tilde \alpha^{-1} \tilde \Delta \tilde \alpha - 2 \Lie_{\tilde \nu} \tilde K + \frac{n+1}{n} \tilde K^2 =
  \frac{4n[-\tr^4 + (n-5)\tr^2 + n]}{(\tr^2+1)^3}.
\end{equation}


\subsection{Wave equation}
\label{s:waveeq}

The scalar field $\Phi$ is supposed to obey the NLW 
\begin{equation} \label{e:nlw2}
  \Box \Phi = \mu |\Phi|^{p-1} \Phi,
\end{equation}
where $\Box$ is the d'Alembertian of the Minkowski metric $\eta$.

Defining a conformally rescaled scalar field 
\begin{equation}  \label{e:Phitildedef}
  \tilde \Phi = \Omega^{(1-n)/2} \Phi
\end{equation}
and using the conformal transformation properties of the d'Alembertian and the Ricci scalar, we have
\begin{equation} \label{e:confNLW}
  \fl \tilde \Box \tilde \Phi - \frac{n-1}{4n} \tilde R \tilde \Phi  = \Omega^{-(n+3)/2} \left( \Box\Phi - \frac{n-1}{4n} R \Phi \right) = \mu \Omega^{[p(n-1)-n-3]/2} |\tilde \Phi|^{p-1} \tilde \Phi,
\end{equation}
where on the right-hand side we have used \eref{e:nlw2} and the fact that $\eta$ is flat, $R = 0$.
This equation is regular at \Scri, where $\Omega=0$, iff $n \geqslant 2$ and
\begin{equation}
  p(n-1)-n-3\geqslant 0 \; \Leftrightarrow \; p \geqslant \frac{n+3}{n-1}
  =: p_\mathrm{conf}.
\end{equation}
For any $n\geqslant 3$ we have $p_\mathrm{conf} < p_\mathrm{crit} $,
where $p_\mathrm{crit}$ was defined in \eref{e:pcrit}, cf. table \ref{t:critexps}.
Thus the conformal approach is capable of treating both energy-subcritical and supercritical NLWs.

\begin{table} 
  \[
  \begin{array}{c||c|c|c|c|c|c|c}
    n & 2 & 3 & 4 & 5 & 6 & 7 & \cdots\\\hhline{=||=|=|=|=|=|=|=}
    p_\mathrm{conf} & 5 & 3 & 7/3 & 2 & 9/5 & 10/6 & \cdots\\\hline
    p_\mathrm{crit} & \infty & 5 & 3 & 7/3 & 2 & 9/5 & \cdots
  \end{array}
\]
  \caption{\label{t:critexps}
    The smallest allowed exponent $p_\mathrm{conf}$ of the nonlinearity for the conformal approach to be applicable, as compared with the energy-critical exponent $p_\mathrm{crit}$, depending on the number of space dimensions $n$.
  For $n=2$ the NLW is energy-subcritical for any $p>1$.}
\end{table}

Performing an $(n+1)$-decomposition of the conformal d'Alembertian, \eref{e:confNLW} becomes
\begin{equation} \label{e:wave2ndorder}
  \fl -\Lie_{\tilde \nu}^2 \tilde \Phi + \tilde \Delta \tilde \Phi + \tilde K \Lie_{\tilde \nu} \tilde \Phi + \tilde \alpha^{-1} \tilde \alpha^{,i} \tilde \Phi_{,i} - \frac{n-1}{4n} \tilde R \tilde \Phi = \mu \Omega^{[p(n-1)-n-3]/2} |\tilde \Phi|^{p-1} \tilde \Phi,
\end{equation} 
where spatial indices $i$ are raised and lowered using the flat conformal metric $\tilde \gamma$.
We can write this equation in first-order in time form by introducing an auxiliary field $\tilde \Pi$,
\begin{eqnarray}
  \fl \Lie_{\tilde \nu} \tilde \Phi &=:& \tilde \Pi,\\
  \fl \Lie_{\tilde \nu} \tilde \Pi &=& \tilde \Delta \tilde \Phi + \tilde K \tilde \Pi + \tilde \alpha^{-1} \tilde \alpha^{,i} \tilde \Phi_{,i} 
  - \frac{n-1}{4n} \tilde R \tilde \Phi - \mu \Omega^{[p(n-1)-n-3]/2} |\tilde \Phi|^{p-1} \tilde \Phi.
\end{eqnarray}
Writing out the Lie derivative $\Lie_{\tilde \nu} = \tilde \alpha^{-1} (\partial_{\tilde t} - \Lie_{\tilde \beta})$, we arrive at
\begin{eqnarray}
  \label{e:scalarwave1}
  \fl \tilde \Phi_{,\tilde t} = \tilde \beta^{\tr} \tilde \Phi_{,\tr} + \tilde \alpha \tilde \Pi, \\
  \fl \tilde \Pi_{,\tilde t} = \tr^{1-n}\left[ \tr^{n-1} \left(\tilde \beta^{\tr} \tilde \Pi + \tilde \alpha \tilde \Phi_{,\tr} \right) \right]_{,\tr} +\tilde \alpha \tr^{-2} \mathring{\Delta}^{(n-2)} \tilde \Phi 
  - \frac{n-1}{4n} \tilde \alpha \tilde R \tilde \Phi \nonumber\\ \fl \qquad 
  - \mu \tilde \alpha \Omega^{[p(n-1)-n-3]/2} |\tilde \Phi|^{p-1} \tilde \Phi.
  \label{e:scalarwave2}
\end{eqnarray} 

We have double-checked using computer algebra that \eref{e:scalarwave1}--\eref{e:scalarwave2} are equivalent to \eref{e:nlw2} if the coordinate transformation given by \eref{e:hyperboloid} and \eref{e:rofr_} is applied directly.

Equations \eref{e:scalarwave1}--\eref{e:scalarwave2} are perfectly regular at future null infinity $\tr=1$, where $\Omega=0$, provided that $p\geqslant p_\mathrm{conf}$ (or $\mu=0$).
No boundary conditions are needed there because all characteristics point towards the exterior of the domain.

It should be noted that we have written the radial principal part of \eref{e:scalarwave2} in conservative form. 
The finite differencing (section \ref{s:spatialdiscr}) will be applied in precisely this order.
This was found to be essential for the numerical evolutions to be stable.


\subsection{Energy balance}
\label{s:energy_theory}

On standard slices of Minkowski spacetime approaching spacelike infinity, the energy \eref{e:stdenergy} is conserved.
In contrast, the energy computed on hyperboloidal slices decreases during the evolution.
This is due to the fact that hyperboloidal slices intersect \Scri, so that outgoing waves pass through this conformal boundary, carrying away energy.
We will compute this energy flux explicitly.
A similar energy balance was derived in \cite{Racz2011} in the context of linear scalar fields in Kerr spacetime, where an additional inner boundary arises at the black hole horizon.

The scalar field $\Phi$ is associated with an energy-momentum tensor 
\begin{equation}
  T_{ab} = \nabla_a \Phi \nabla_b \Phi - \half \eta_{ab} \nabla_c \Phi \nabla^c \Phi - \frac{\mu}{p+1} |\Phi|^{p+1} \eta_{ab},
\end{equation}
where $\nabla$ is the covariant derivative of the Minkowski metric $\eta$.
The vanishing of its divergence, $\nabla^b T_{ab}=0$, is equivalent with the NLW \eref{e:nlw2}.
The Killing vector $k := \partial/\partial t = \partial/\partial \tilde t$ gives rise to the conserved current
\begin{equation}
  E^a = T^a{}_b k^b = T^a{}_t
\end{equation} 
satisfying
\begin{equation} \label{e:divE}
  \nabla_a E^a = 0. 
\end{equation}
Let $\Sigma_{\tilde t}$ denote the slice of constant hyperboloidal time $\tilde t$ with future-directed unit normal $\nu$ (normalised w.r.t. $\eta$) and induced Riemannian metric $\gamma_{ab} = \eta_{ab} + \nu_a \nu_b$.
Furthermore let $s$ denote the outward-directed normal to the timelike surface $r=\const$ (again normalised w.r.t. $\eta$) and $h_{ab} = \eta_{ab} - s_a s_b$ the induced pseudo-Riemannian metric on such a surface.
Applying Gauss' theorem to \eref{e:divE}, the energy balance reads
\begin{eqnarray} \label{e:energybalance}
  \fl E(\tilde t_2) - E(\tilde t_1) := \int_{\Sigma_{\tilde t_2}}  \nu_a E^a \sqrt{\det \gamma} \,\rmd r \, \rmd^{n-1} \theta - \int_{\Sigma_{\tilde t_1}} \nu_a E^a \sqrt{\det \gamma} \, \rmd r \, \rmd^{n-1}\theta \nonumber\\ 
  \fl \qquad = \int_{\tilde t_1}^{\tilde t_2} \int_{S^{n-1}} s_a E^a \sqrt{\det |h|} \, \rmd^{n-1} {\bf \theta} \, \rmd \tilde t \, \Big\vert_{\mathrsfs{I}^+} =: F(\tilde t_1, \tilde t_2),
\end{eqnarray}
where $\theta$ stands collectively for the $n-1$ angular coordinates on $S^{n-1}$.

Expressing $\Phi$ in terms of the conformally rescaled scalar field $\tilde\Phi$ defined in \ref{e:Phitildedef}, we obtain
for the energy
\begin{eqnarray} \label{e:energy}
  \fl E(\tilde t) = \frac{C}{4n}\int_0^1 \tr^{n-1} \rmd \tr \int_{S^{n-1}} dS^{(n-1)}  
  \left\{ (1+\tr^2) \left[ \tilde \Pi^2 + \tilde \Phi_{,\tr}^2 + \frac{1}{\tr^2} \Vert \mathring{\nabla}^{(n-1)} \tilde \Phi \Vert^2 \right.\right.
  \nonumber\\ 
  \left. \fl\qquad + 2\mu \Omega^{[p(n-1)-n-3]/2} \frac{1}{p+1} |\tilde \Phi|^{p+1} \right] 
  \nonumber\\ \fl\qquad
  \left. + 2(n-1) \frac{\tr^2-1}{\tr^2+1} \, \tr \tilde \Phi_{,\tr}\tilde \Phi - 4 r \tilde \Phi_{,\tr} \tilde \Pi + (n-1)^2 \frac{\tr^2}{\tr^2+1} \tilde \Phi^2 \right\}.
\end{eqnarray}
Here $\mathring{\nabla}^{(n-1)}$ refers to the gradient on $S^{n-1}$ and 
$dS^{(n-1)} = \sqrt{\det \sigma^{(n-1)}} \rmd^{n-1} \theta$ to its area element, e.g. for $n=3$
\begin{equation}
  \Vert \mathring{\nabla}^{(2)} \Phi \Vert^2 = \Phi_{,\theta}^2 + \frac{1}{\sin^2 \theta} \, \Phi_{,\varphi}^2, \qquad
    dS^{(2)} = \sin\theta \, \rmd\theta \,\rmd\varphi. 
\end{equation}
In the first line of \eref{e:energy}, we recognize the usual energy expression for the linear scalar field.
The second line of \eref{e:energy} contains the contribution from the nonlinear potential, and the terms in the third line arise from the conformal rescaling.

The flux integral takes the form
\begin{equation} \label{e:flux}
  F(\tilde t_1, \tilde t_2) = -\frac{C^2}{n^2} \int_{\tilde t_1}^{\tilde t_2} \rmd \tilde t\int_{S^{(n-1)}} dS^{(n-1)}  (\tilde \Phi_{,\tr} - \tilde \Pi)^2 \, \Big \vert_{\tr=1}.
\end{equation}
It is manifestly negative, reflecting the fact that the waves carry away energy as they pass through \Scri.


\section{Numerical methods} 
\label{s:nummeth}

In this section we describe the numerical methods we use in order to solve the NLW \eref{e:scalarwave1}--\eref{e:scalarwave2}.
We begin with the spatial discretisation in section \ref{s:spatialdiscr}.
Next we explain how the equation is integrated forward in time in section \ref{s:timestepping}, and we address how to eliminate high-frequency instabilities.
Finally we provide details on how spatial integrals needed for various diagnostics are computed in section \ref{s:integrals}.
The code has been written in Python using libraries such as \texttt{NumPy} \cite{NumPy} and \texttt{SciPy} \cite{SciPy}.


\subsection{Spatial discretisation}
\label{s:spatialdiscr}

We use a hybrid method for the spatial discretisation, which combines a finite-difference method in the radial direction with a pseudo-spectral method in the angular directions.
The latter has the advantage that it is easier to build in parity conditions by an appropriate choice of expansion functions, as will be explained below.
Furthermore, if we introduce a logically rectangular grid with equidistant grid points in the standard spherical coordinates $\theta$ and $\varphi$ on the two-sphere, the distance between neighbouring $\varphi$-grid points becomes very small near the poles. 
This leads to a severe restriction on the time step due to the Courant-Friedrichs-Lewy (CFL) condition if a finite-difference method is used.
With a pseudo-spectral method, it is straightforward to remove an increasing number of higher (more oscillatory) $\varphi$-modes as the axis is approached (section \ref{s:timestepping}), thus allowing for larger time steps.

We take the spatial dimension to be $n=3$ without symmetries first. 
Higher dimensions with $\mathrm{SO}(n-1)$ symmetry will be treated as a special case.

We introduce equidistant grid points in each dimension,
\begin{eqnarray}
  \tr_i = (i + \half) \htr, \quad\; 
  i=0,1,\ldots, \Ntr - 1, \quad 
  h_{\tr} = 1/(\Ntr - \half),\\
  \theta_j = (j + \half) h_\theta, \quad 
  j=0,1,\ldots, N_\theta - 1, \quad 
  h_\theta = \pi/N_\theta,\\
  \varphi_k = k h_\varphi, \qquad\;\; 
  k = 0, 1, \ldots, N_\varphi - 1, \quad 
  h_\varphi = 2\pi/N_\varphi.
\end{eqnarray}
The radial grid is staggered at the origin, whereas the outermost grid point coincides with \Scri ($\tr=1)$.
The grid in the elevation $\theta$ is staggered about the axis.
The reason for this choice is that there are terms in the wave equation that are formally singular at the origin ($\tr=0$) and on the axis ($\sin\theta=0$), which would be difficult to evaluate if there were grid points at the origin and on the axis.

Fourth-order finite differences are used in order to approximate the radial derivatives. 
In the interior a centred stencil is used.
Setting $u_i := u(\tr_i)$ for a smooth function $u$ and leaving out its angular dependence for simplicity, 

\begin{equation}
  u'_i \approx (u_{i-2} - 8 u_{i-1} + 8 u_{i+1} - u_{i+2})/(12 \htr)
\end{equation}
for $i=0,1,\ldots,\Ntr-3$, where $\approx$ stands for equality up to a truncation error of the 
order $O(\htr^4)$.
Evaluating the above stencils at $i=0$ and $i=1$ requires values of $u$ at the ghost points 
$\tr_{-1}=-h_{\tr}/2$ and $\tr_{-2}=-3h_{\tr}/2$, which we copy from interior points according to 
\begin{equation}
  \label{e:symmetryorigin}
  u(-\tr, \theta, \varphi) = u(\tr, \pi-\theta, \pi+\varphi),
\end{equation}
satisfied by any smooth function $u$.
We require $N_\varphi$ to be even, which ensures that if there is a grid point at $\varphi$, there is also a grid point at $\pi+\varphi$.

Near the outer boundary at $\tr=1$, one-sided differences are used,
\begin{eqnarray}
  \fl u'_{\Ntr-2} \approx (-u_{\Ntr-5} + 6 u_{\Ntr-4} - 18 u_{\Ntr-3} + 10 u_{\Ntr-2} + 3 u_{\Ntr-1})/(12 \htr),\\
  \fl u'_{\Ntr-1} \approx (3 u_{\Ntr-5} - 16 u_{\Ntr-4} + 36 u_{\Ntr-3} - 48 u_{\Ntr-2} + 25 u_{\Ntr-1} )/(12 \htr),
\end{eqnarray}
the truncation error again being $O(\htr^4)$.

Our pseudo-spectral method is based on Fourier expansions in both $\theta$ and $\varphi$. 
This allows us to employ fast Fourier transform (FFT) techniques \cite{Boyd1978}, which would not be available if we expanded the fields in spherical harmonics.

Leaving out its radial dependence for simplicity, a smooth function $u$ is expanded as
\begin{equation} \label{e:spectralexpansion}
  \fl u(\theta,\varphi) \approx \sum_{l=0}^{N_\theta-1} \left( \cos(l\theta) \sum_{m=0\atop{m\,\mathrm{even}}}^{N_\varphi/2-1} a_{lm} \, \rme^{\rmi m \varphi} +  \sin(l\theta) \sum_{m=1\atop{m\,\mathrm{odd}}}^{N_\varphi/2-1} a_{lm} \, \rme^{\rmi m \varphi} \right).
\end{equation}
Any smooth function $u$ on the two-sphere must obey 
\begin{equation}
  u(-\theta, \varphi) = u(\theta, \pi + \varphi), \qquad
  u(\pi + \theta, \varphi) = u(\pi - \theta, \pi + \varphi),
\end{equation}
and the way we expand the even and odd $m$-modes separately in \eref{e:spectralexpansion} enforces this.
For $u$ to be real, the coefficients $a_{lm}$ in \eref{e:spectralexpansion} must be complex, and we may use 
\begin{equation}
  \Re(a_{lm} \rme^{\rmi m \varphi}) = \Re \, a_{lm} \cos(m\varphi) - \Im \, a_{lm} \sin(m\varphi).
\end{equation}

In order to compute the expansion coefficients $a_{lm}$ from the point values $u_{jk} := u(\theta_j, \varphi_k)$, we first perform a real FFT w.r.t. $\varphi$ using \texttt{scipy.fft.rfft}.
Next, we perform a discrete cosine transform w.r.t. $\theta$ for the even-$m$ modes (\texttt{scipy.fft.dct}) and a discrete sine transform for the odd-$m$ modes (\texttt{scipy.fft.dst}).
The (default) Type II versions of the discrete cosine and sine transforms implemented in \texttt{scipy} are appropriate for our staggered $\theta$ grid.

Derivatives of the expansion can now be computed analytically by simply differentiating the basis functions $\cos(l\theta)$, $\sin(l\theta)$ and $\rme^{\rmi m \varphi}$.
Finally in order to obtain the values of the derivatives at the grid points $(\theta_j, \varphi_k)$, we apply the inverse discrete cosine and sine transforms w.r.t. $\theta$, followed by an inverse real FFT w.r.t. $\varphi$.

In the code the main object is the array of values of the unknowns $\tilde \Phi$ and $\tilde \Pi$ at the grid points $(\tr_i, \theta_j, \varphi_k)$. 
Terms in the evolution equations, especially nonlinear terms, are evaluated and added pointwise.

We remark that a purely spectral (rather than pseudo-spectral collocation) method such as the one used in \cite{Racz2011} cannot be employed in our case because of the nonlinear term in the wave equation we consider.

\paragraph{$n$ dimensions with $\mathrm{SO}(n-1)$ symmetry.}
In this case there is no $\varphi$-dependence and the expansion \eref{e:spectralexpansion} is replaced with
\begin{equation} \label{e:spectralexpansionreduced}
  u(\theta) \approx \sum_{l=0}^{N_\theta-1} a_l \cos(l\theta) 
\end{equation}
with real coefficients $a_l$.
This satisfies the symmetry conditions 
\begin{equation}
  u(-\theta) = u(\theta), \qquad
  u(\pi + \theta) = u(\pi - \theta).
\end{equation}


\subsection{Time stepping and filtering}
\label{s:timestepping}

We employ the method of lines, whereby the evolution equations \eref{e:scalarwave1}--\eref{e:scalarwave2} are first discretised in space as described in the previous subsection, leading to a system of ordinary differential equations (ODEs), one for each grid point.
These ODEs are then integrated forward in time using the standard fourth-order Runge-Kutta method.

When using standard finite differences in the method of lines, the resulting system is generally not numerically stable \cite{Kreiss1973,Gustafsson2013}.
There are high-frequency modes that can grow without bound.
In order to eliminate these, we apply fifth-order Kreiss-Oliger dissipation \cite{Kreiss1973} in the radial direction, whereby the term 
\begin{eqnarray}
  \fl (Q u)_i = \frac{\varepsilon}{64 \htr} (u_{i-3} - 6 u_{i-2} + 15 u_{i-1} - 20 u_i + 15 u_{i+1} - 6 u_{i+2} + u_{i+3}) \nonumber \\ 
  \fl \qquad \; = \frac{\varepsilon}{64} \htr^5 \left(\frac{\rmd^6 u}{\rmd\tilde r^6}\right)_i + \Or(\htr^7) \label{e:diss} 
\end{eqnarray}
for $i=0,1,\ldots,\Ntr-4$ is added to the right-hand side of the discretised evolution equations for both $\tilde\Phi$ and $\tilde\Pi$.
Ghost points near the origin are filled according to \eref{e:symmetryorigin}.
No dissipation is added at the outermost three radial grid points.
We typically take the parameter $\varepsilon = 0.2$.
It should be noted that the artificial dissipation term \eref{e:diss} is of higher order in $h_{\tr}$ than the truncation error of the finite-difference method.

In a pseudo-spectral method, nonlinear terms such as the one in our wave equation can cause high frequencies to be represented incorrectly due to aliasing \cite{Boyd2000}.
We address this by applying spectral filtering according to the Orszag 2/3 rule \cite{Patterson1971,Boyd2000}.
At the end of each of the four substeps of the fourth-order Runge-Kutta method, we transform the solution to Fourier space w.r.t. $\theta$ and $\varphi$, i.e., we compute the coefficients $a_{lm}$ in \eref{e:spectralexpansion} (or $a_l$ in \eref{e:spectralexpansionreduced}).
We then set the top third w.r.t. $l$ of these coefficients to zero. 
Finally we transform back to the point values.
We remark that while strictly speaking, the 2/3 rule can only be shown mathematically to eliminate aliasing errors from quadratic nonlinearities, it is widely used in practice for various nonlinear partial differential equations, and we have not observed any high-frequency instabilities when applying it.
 
In order to combat the clustering of grid points near the poles of the two-sphere, we filter out an increasing number of the highest Fourier modes w.r.t. $\varphi$ as $\theta$ approaches $0$ or $\pi$.
Specifically, we remove the proportion $1 - \sin\theta$ of all $\varphi$-Fourier modes (beginning at the highest frequencies) \cite{Fornberg1998}.

The time step is restricted by the smallest distance between neighbouring grid points, and this occurs between neighbouring $\theta$-grid points at the smallest radius:
$\Delta x_{\min} = \tr_0 h_\theta$, where $\tr_0 = \htr/2$.
The time step is then taken to be $\Delta \tilde t = \lambda \Delta x_{\min}$,
where $0<\lambda<1$ according to the CFL condition.
We typically choose $\lambda = 0.8$.


\subsection{Spatial integration}
\label{s:integrals}

For expressions such as the energy and flux (section \ref{s:energy}) we will need to compute integrals on a $\tilde t=\const$ slice and on its spherical boundary.
The radial integration is performed using Simpson's rule (taking $\Ntr$ to be even), 
\begin{eqnarray}
  \fl \int_{0}^1 u(\tr) \rmd \tr \approx \frac{\htr}{3} \left(\frac{27}{8} u_0 + \frac{17}{8} u_1 + 4 u_2 + 2 u_3 + 4 u_4 + 2 u_5 + \ldots 
  \right. \nonumber\\ \left. + 4 u_{\Ntr-4} + 2 u_{\Ntr-3} + 4 u_{\Ntr-2} + u_{\Ntr-1}\right),
\end{eqnarray}
where the error is of the same order as the truncation error of the finite-difference scheme, $\Or(\htr^4)$.
The modified coefficients of $u_0$ and $u_1$ are due to the fact that the grid is staggered at the origin, with a first grid point at $\tr_0 = \htr/2$, and $u$ is assumed to be an even function of $\tr$ here.

The spectral expansion on the right-hand side of \eref{e:spectralexpansion} can easily be integrated exactly:
\begin{equation}
   \int_{S^{2}} u \, dS^{(2)} =
   \int_0^{2\pi} \int_0^\pi u(\theta,\varphi) \, \sin\theta \, \rmd \theta \, \rmd \varphi \approx 2\pi \sum_{l=0\atop{l\,\mathrm{even}}}^{N_\theta - 1}\frac{2a_{l0}}{1-l^2}.
\end{equation} 

\paragraph{$n$ dimensions with $\mathrm{SO}(n-1)$ symmetry.}
In this case \eref{e:spectralexpansionreduced} integrates to
\begin{equation} \fl
  \int_{S^{n-1}} u \, dS^{(n-1)} = \int_{S^{n-2}} \int_0^\pi u(\theta) \sin^{n-2}\theta \, \rmd\theta \, dS^{(n-2)} \approx A^{(n-2)} \sum_{l=0}^{N_\theta-1} a_l \, I_{n-2, l},
\end{equation} 
where 
\begin{equation} \label{e:areaS}
  A^{(n-2)} = \frac{2\pi^{(n-1)/2}}{\Gamma(\tfrac{n-1}{2})}
\end{equation}
is the area of $S^{n-2}$, and the integrals 
\begin{equation}
  I_{n-2, l} := \int_0^\pi \sin^{n-2} \theta \, \cos(l \theta) \, \rmd\theta
\end{equation}
are computed numerically to high precision using \texttt{scipy.integrate.quad} in the initialisation phase of the code.
Alternatively, one may use\footnote{
  The author is grateful to an anonymous referee for pointing out this formula.}
\begin{equation}
  I_{n-2, l} = \frac{2^{2 - n} \, \pi \cos\left(\frac{l \pi}{2}\right) \, \Gamma(n-1)}{\Gamma\left(\frac{n-l}{2}\right) \, \Gamma\left(\frac{n+l}{2}\right)}.
\end{equation}


\section{Results}
\label{s:results}

In this section we present our numerical results.
We first perform a convergence test against exact solutions to the linear wave equation in section \ref{s:convtest}.
Next we check that the energy balance the scalar field must obey on hyperboloidal slices is satisfied in our numerical evolutions (section \ref{s:energy}).
Finally in section \ref{s:tails} we investigate the late-time power-law decay (tails) of the solutions.

Throughout we take the value of the constant $C$ determining the mean curvature of the hyperboloidal slices in \eref{e:hyperboloid} to be $C=0.5$.

 
\subsection{Convergence test against an exact linear solution} 
\label{s:convtest}

Given the complexity of the algebraic operations involved in deriving the evolution equations, it is very useful to have an exact solution at hand in order to check that the implementation is correct, and that the numerical method converges as expected.

In \ref{s:exact} we construct exact solutions to the \emph{linear} wave equation. 
They are based on an expansion in spherical harmonics and are given in terms of a mode function, which we take to be 
\begin{equation}
  F(x) = A \, x \, \exp \left[-\frac{1}{2} \left( \frac{x}{\sigma}\right)^2\right].
\end{equation}
We choose $A=1$ and $\sigma=1$.

For the convergence test in the case of $n=3$ without symmetries, we consider a superposition of two modes with spherical harmonic indices $(l,m)=(1,1)$ and $(2,-2)$. 
In the case of $n=5$ with SO$(4)$ symmetry, we take two modes with $l=1$ and $2$.

The initial data are set according to $\tilde \Phi_0 := \tilde \Phi(\tilde t=0) = \tilde \Phi_\mathrm{exact} (\tilde t = \tilde t_0)$ and similarly for $\tilde \Pi$, where we choose $\tilde t_0 = -15$ so that the solution is initially ingoing and centred roughly in the middle of the conformal radial interval $[0,1]$.

\begin{figure} 
  \begin{subfigure}[b]{0.5\textwidth}
    \includegraphics[width=\textwidth]{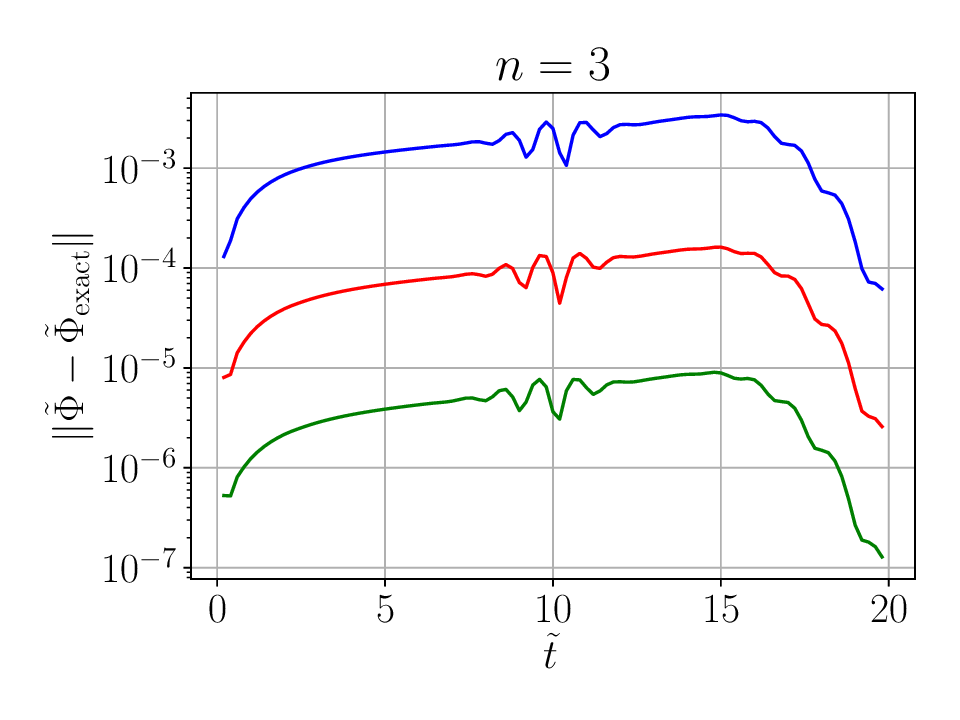}
    \caption{}
  \end{subfigure}
  \begin{subfigure}[b]{0.5\textwidth}
    \includegraphics[width=\textwidth]{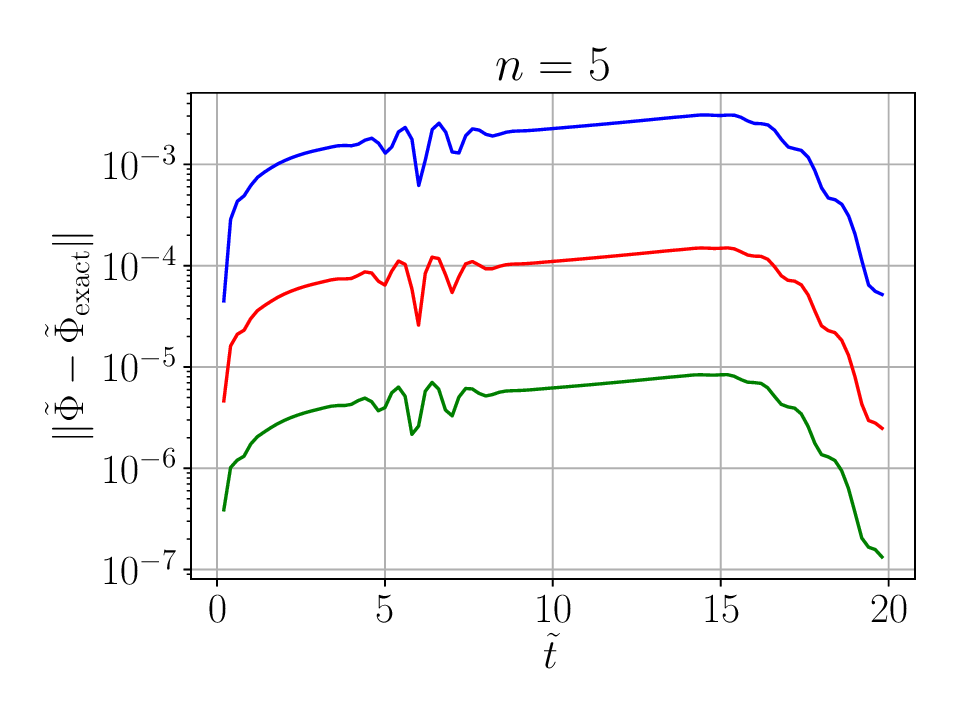}
    \caption{}
  \end{subfigure}
  \caption{\label{f:convtest} 
    Convergence test against an exact linear solution. 
    (a) $n=3$ dimensions without symmetries, (b) $n=5$ with SO$(4)$ symmetry.
    Shown is the $L^2$ norm of the error of $\tilde \Phi$ w.r.t.~the exact linear solution as a function of time for three different radial resolutions, from top to bottom: $\Ntr=250$ (blue), $\Ntr=500$ (red) and $\Ntr=1000$ (green).
  }
\end{figure}

We now evolve these initial data and compare with the exact solution.
In figure \ref{f:convtest} we plot the $L^2$ norm of the error
\begin{equation}
  \Vert \tilde \Phi - \tilde \Phi_\mathrm{exact} \Vert_{L^2} = \left[ \int_0^1 \tr^{n-1} \rmd \tr \int_{S^{n-1}} (\tilde \Phi - \tilde \Phi_\mathrm{exact})^2 \, dS^{(n-1)} \right]^{1/2}
\end{equation}
as a function of time for three different radial resolutions at fixed angular resolution.
We observe that doubling the resolution decreases the error nearly by a factor of $2^4 = 16$, as expected for a fourth-order accurate finite-difference method.

We remark that at the chosen angular resolution, the spherical harmonics are represented exactly by our pseudo-spectral method, so there is no discretisation error in the angular dimensions for these linear evolutions.

No (dynamical) exact solutions are known in the nonlinear case. 
However, in figure \ref{f:energyviolation_conv} in section \ref{s:energy} and figure \ref{f:lpiconv} in section \ref{s:tails} we will compare nonlinear evolutions at different numerical resolutions.


\subsection{Energy balance}
\label{s:energy}

In this section we verify numerically the energy balance derived in section \ref{s:energy_theory} that the scalar field obeys on our hyperboloidal foliation.

We consider two cases: $n=3,\, p=5$ (critical) and $n=5,\, p=3$ (supercritical). 
In both cases we impose an SO$(n-1)$ symmetry. 
The initial data are taken to be momentarily static with
\begin{equation} \label{e:staticID1}
  \tilde \Phi_0 = A \exp \left[-\left( \frac{\tr-\tr_0}{\sigma}\right)^2 \right] \, Y_l(\theta) ,
\end{equation}
and $\tilde\Phi_{,t} = 0$ at $t=0$ implies
\begin{equation} \label{e:staticID2}
  \tilde \Pi_0 = \frac{2\tr}{1+\tr^2} \tilde \Phi_{0,\tr} .
\end{equation}
We superpose two modes with $l=2$ and $l=3$ and choose the parameters
$\tr_0 = 0.3$ and $\sigma=0.07$.
The amplitude for both modes is taken to be $A=12$ for $n=3,\, p=5$ and $A=200$ for $n=5,\, p=3$, which is smaller than but close to the critical amplitude beyond which blow-up occurs in the focusing case ($\mu = -1$).
The numerical resolution is taken to be $\Ntr=4000$, $N_\theta=24$.

\begin{figure} 
  \begin{subfigure}[b]{0.5\textwidth}
    \includegraphics[width=\textwidth]{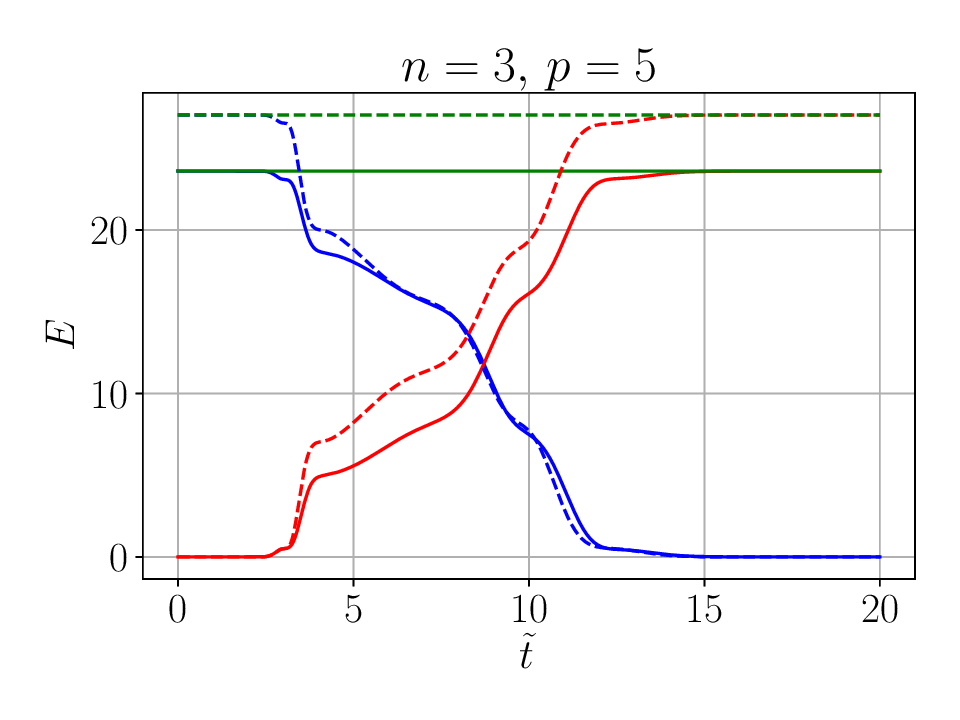}
    \caption{}
  \end{subfigure}
  \begin{subfigure}[b]{0.5\textwidth}
    \includegraphics[width=\textwidth]{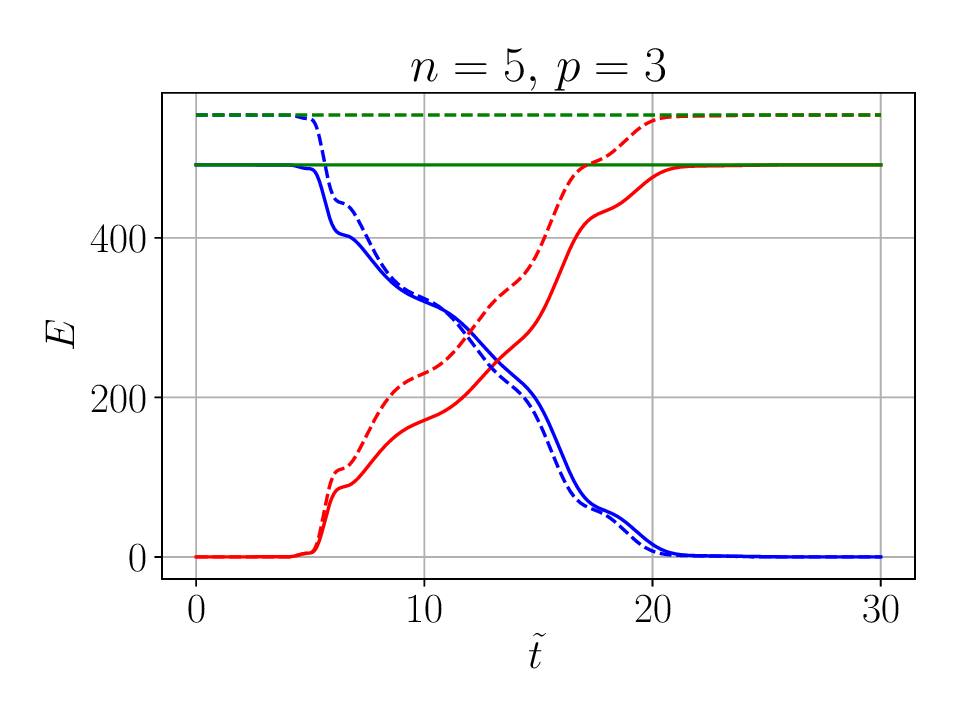}
    \caption{}
  \end{subfigure}
  \caption{\label{f:energybalance}  
    Energy balance for the nonlinear wave equation. (a) $n=3,\, p=5$, (b) $n=5,\, p=3$, both with SO$(n-1)$ symmetry.
    Shown is the energy $E(\tilde t)$ (blue, monotonically decreasing), the integrated flux $-F(0,\tilde t)$ (red, monotonically increasing) and their sum (green, nearly constant).
    Solid lines correspond to an evolution with a focusing nonlinearity ($\mu=-1$), dashed lines to a defocusing nonlinearity ($\mu=1$).
  }
\end{figure}

Figure \ref{f:energybalance} demonstrates that the energy balance (cf. equation \eref{e:energybalance})
\begin{equation}
  E(\tilde t) - F(0,\tilde t) = E(0)
\end{equation}
is well satisfied numerically, i.e., the sum of the numerically computed terms on the left-hand side is indeed nearly constant.
For an evolution with a defocusing nonlinearity ($\mu=1$), the initial energy is higher than for a focusing nonlinearity ($\mu=-1$) due to the different sign of the potential energy term in \eref{e:energy}.
 
\begin{figure} 
  \begin{subfigure}[b]{0.5\textwidth}
    \includegraphics[width=\textwidth]{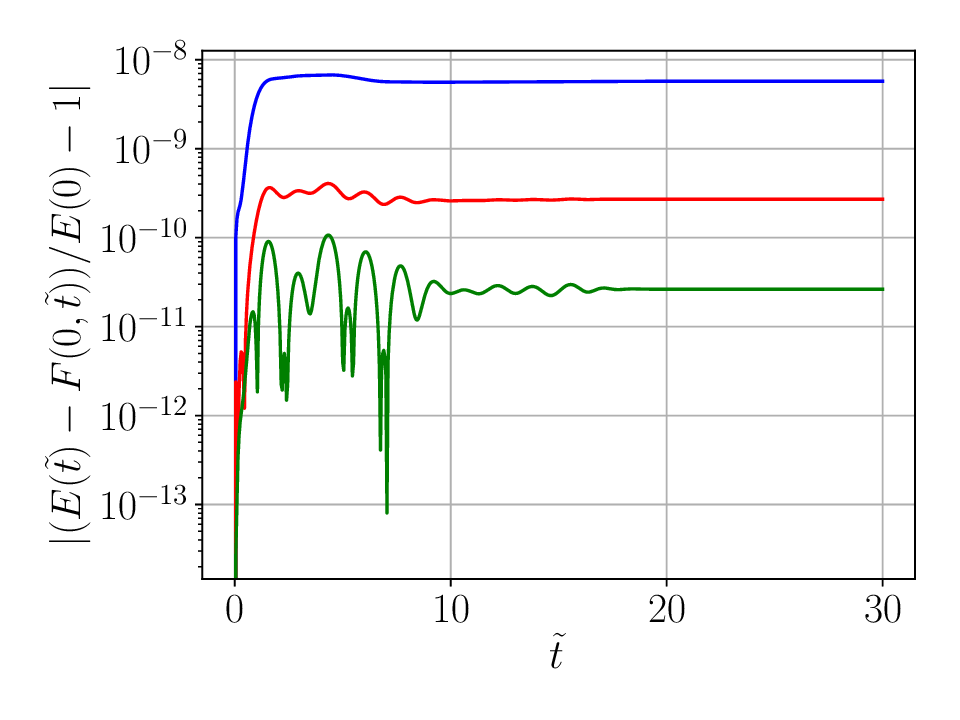}
    \caption{}
  \end{subfigure}
  \begin{subfigure}[b]{0.5\textwidth}
    \includegraphics[width=\textwidth]{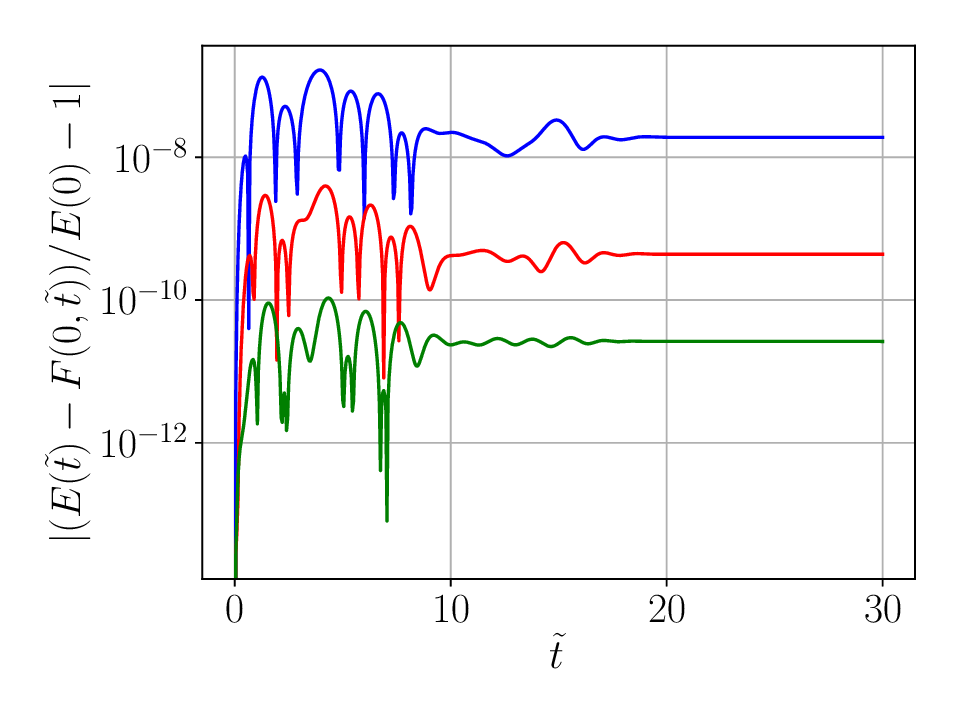}
    \caption{}
  \end{subfigure}
  \caption{\label{f:energyviolation_conv}  
  Convergence test for the relative error in the energy balance of a defocusing $n=5$, $p=3$ NLW evolution.
  (a) Varying radial resolution $\Ntr=1000, 2000, 4000$ (from top to bottom) at fixed angular resolution $N_\theta=24$, 
  (b) varying angular resolution $N_\theta=16, 20, 24$ (from top to bottom) at fixed $\Ntr=4000$.
  }
\end{figure}

In figure \ref{f:energyviolation_conv} we plot the relative error in the energy balance,
\begin{equation}
  \frac{E(\tilde t) - F(0,\tilde t) - E(0)}{E(0)},
\end{equation}
as a function of time for the $n=5$, $p=3$ defocusing case.
Nearly fourth-order convergence of this error towards zero can be observed as the radial resolution is increased, and nearly exponential convergence as the angular resolution is increased, as is typical of (pseudo-)spectral methods.

\begin{figure}
  \begin{subfigure}[b]{0.5\textwidth}
    \includegraphics[width=\textwidth]{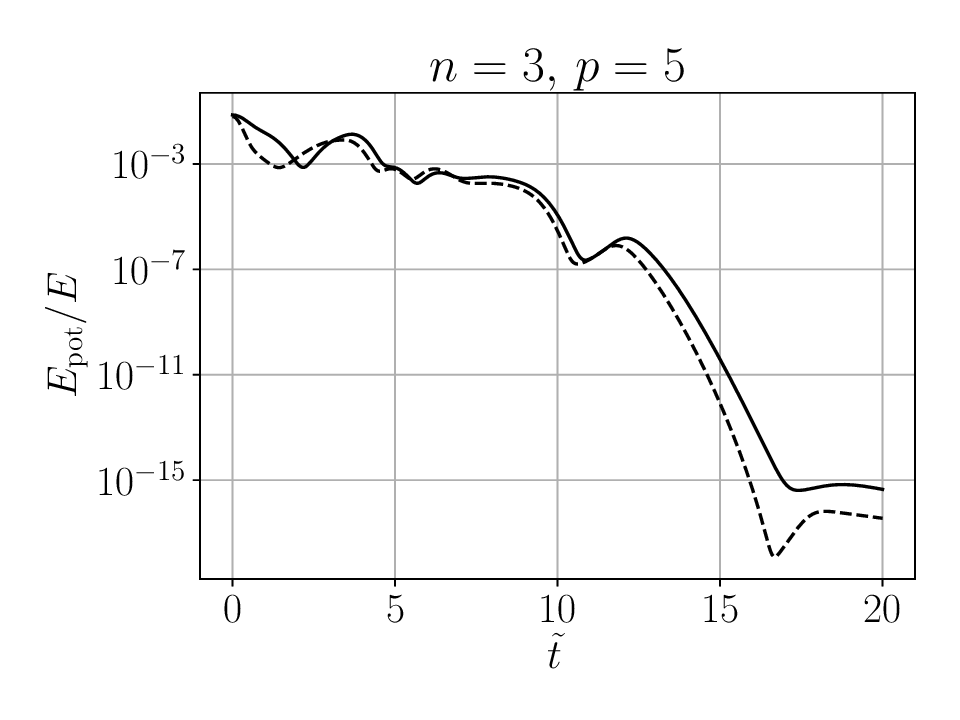}
    \caption{}
  \end{subfigure}
  \begin{subfigure}[b]{0.5\textwidth}
    \includegraphics[width=\textwidth]{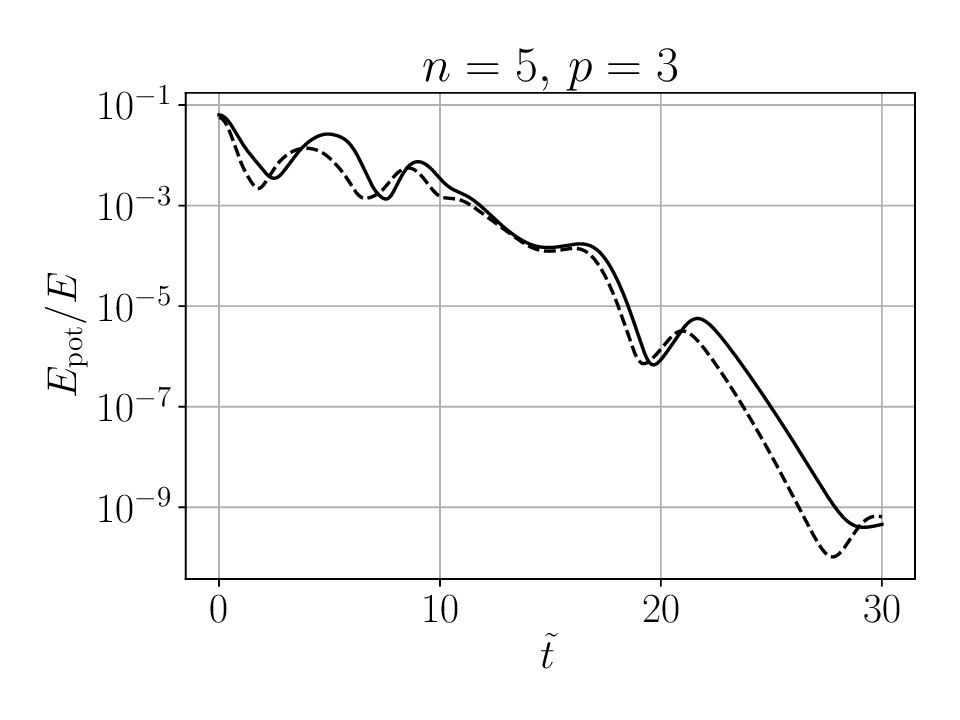}
    \caption{}
  \end{subfigure}
  \caption{\label{f:energyratio} 
    Ratio $E_\mathrm{pot}(\tilde t)/E(\tilde t)$ of the potential energy to the total energy as a function of time for the same evolutions as in figure \ref{f:energybalance}.
    (a) $n=3,\, p=5$, (b) $n=5,\, p=3$, both with SO$(n-1)$ symmetry.
    Solid lines correspond to an evolution with a focusing nonlinearity ($\mu=-1$), dashed lines to a defocusing nonlinearity ($\mu=1$).
  } 
\end{figure}

We also compute the potential energy contribution to \eref{e:energy},
\begin{equation}
  \fl E_\mathrm{pot}(\tilde t) = \frac{C\mu}{2n}\int_0^1 \tr^{n-1} \rmd \tr \int_{S^{n-1}} dS^{(n-1)} (1+\tr^2) \Omega^{[p(n-1)-n-3]/2} \frac{1}{p+1} |\tilde \Phi|^{p+1} ,
\end{equation}
and its ratio to the total energy, $E_\mathrm{pot}(\tilde t)/E(\tilde t)$.
The result is shown in figure \ref{f:energyratio} and indicates that the potential energy becomes negligible at late times.
 
\subsection{Tails} 
\label{s:tails}

In this section we focus on the late-time behaviour of the scalar field, the so-called power-law tails.
We extract the scalar field at various radii and decompose it into spherical harmonics. 
Each mode shows a characteristic power-law decay $\sim \tilde t^{-q}$ at late times.
We investigate how the decay exponent $q$ depends on the power $p$ of the nonlinearity, the spherical harmonic mode indices $(l,m)$ and the extraction radius.

At a fixed extraction radius $\tr_\mathrm{ex}$, we expand the scalar field in spherical harmonics, here in $n=3$ dimensions without symmetries:
\begin{equation}
  \tilde \Phi \vert_{\tr=\tr_\mathrm{ex}} = \sum_{l=0}^\infty \sum_{m=-l}^l \tilde \Phi_{lm} (\tilde t) \, \hat Y_{l m}(\theta, \varphi).
\end{equation}
Here $\hat Y_{l m}$ refers to a real basis of spherical harmonics obtained from the original complex spherical harmonics $Y_{lm}$ via equation \eref{e:real_sph_harm} in \ref{s:exactn3}.
By orthonormality, the mode functions $\tilde \Phi_{lm} (\tilde t)$ can be computed via integration on the sphere using the numerical techniques of section \ref{s:integrals}:
\begin{equation}
  \tilde \Phi_{lm} (\tilde t) = \int_0^{2\pi}\int_0^\pi   \tilde \Phi(\tilde t, \tr, \theta,\varphi) \hat Y_{l m} (\theta, \varphi) \sin\theta \, \rmd\theta\, \rmd\varphi \Big\vert_{\tr=\tr_\mathrm{ex}}.
\end{equation}

For each of the modes we form the local power index
\begin{equation} \label{e:lpidef}
  q_{lm}(\tilde t) := -\frac{\rmd \ln \tilde \Phi_{lm}}{\rmd \ln \tilde t} \Big\vert_{\tr=\tr_\mathrm{ex}} 
  = -\frac{\tilde t \, (\tilde \Phi_{,\tilde t})_{lm}}{\tilde \Phi_{lm}} \Big\vert_{\tr=\tr_\mathrm{ex}}.
\end{equation}
If this approaches a constant $q_{lm}$ as $\tilde t\to\infty$ then the mode decays asymptotically as $\tilde \Phi_{lm} \sim  \tilde t^{-q_{lm}}$.

\begin{figure} 
  \begin{subfigure}[b]{0.5\textwidth}
    \includegraphics[width=\textwidth]{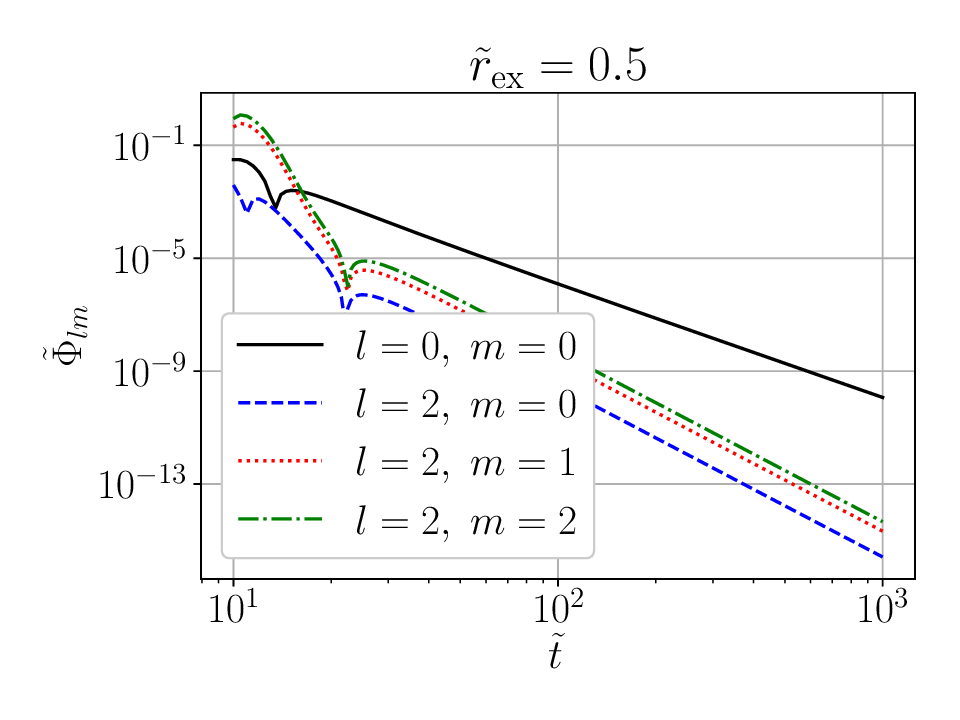}
    \caption{}
  \end{subfigure}
  \begin{subfigure}[b]{0.5\textwidth}
    \includegraphics[width=\textwidth]{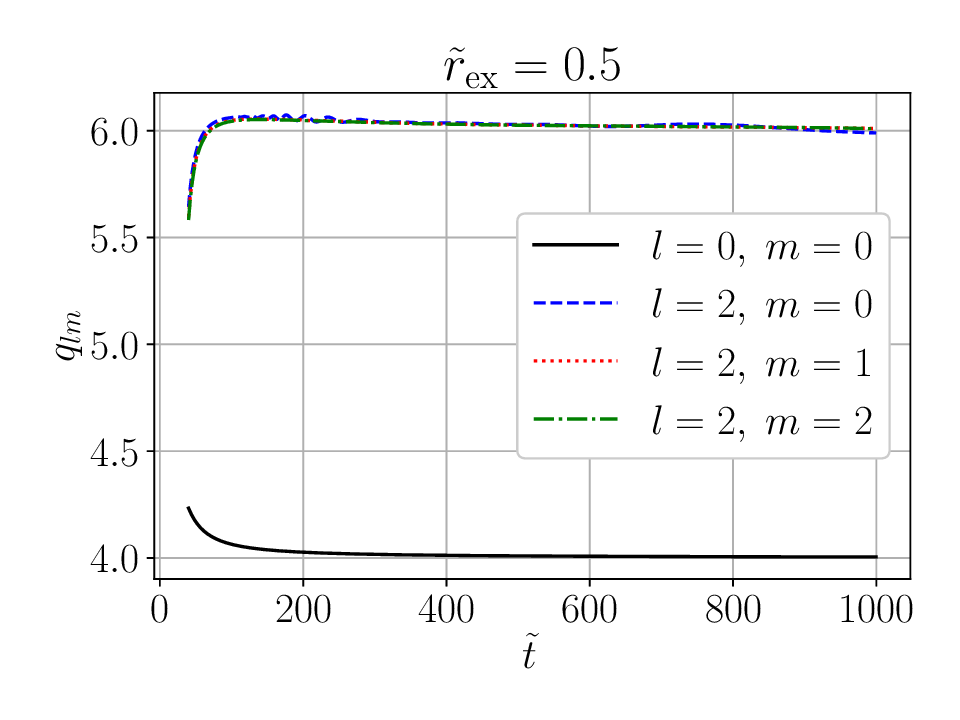}
    \caption{}
  \end{subfigure}
  \caption{\label{f:tails3d} 
    (a) Modes $\tilde\Phi_{lm}$ and (b) local power indices $q_{lm}$ extracted at $\tr_\mathrm{ex}=0.5$ for an evolution in $n=3$ dimensions with a focusing $p=5$ nonlinearity.
    Static initial data containing two modes with $l=2,\, m=1$ (amplitude $A=6$) and $l=2,\,m=2$ (amplitude $A=12$) are chosen. 
    The remaining initial data parameters are $\tr_0 = 0.3$ and $\sigma = 0.07$ for both modes. 
    The numerical resolution is $\Ntr=2000,\, N_\theta=N_\varphi=8$.
  }
\end{figure}

Our first observation is that the asymptotic decay rate appears to be independent of the azimuthal spherical harmonic index $m$ in three spatial dimensions.
In order to illustrate this, figure \ref{f:tails3d} shows an evolution of static initial data \eref{e:staticID1}-\eref{e:staticID2} containing two modes with $(l,m)=(2,1)$ and $(2,2)$. 
Via the nonlinearity ($p=5$, focusing), modes with $(l,m)=(0,0)$ and $(2,0)$ are also excited during the evolution.
The $l=2$ modes all decay at the same rate, $q_{2m}=6$ irrespective of $m$.

In the following we therefore focus on the $m=0$ case only, i.e. we impose axisymmetry for $n=3$, and more generally an SO$(n-1)$ symmetry in $n$ dimensions.

\begin{figure} 
  \begin{subfigure}[b]{0.5\textwidth}
    \includegraphics[width=\textwidth]{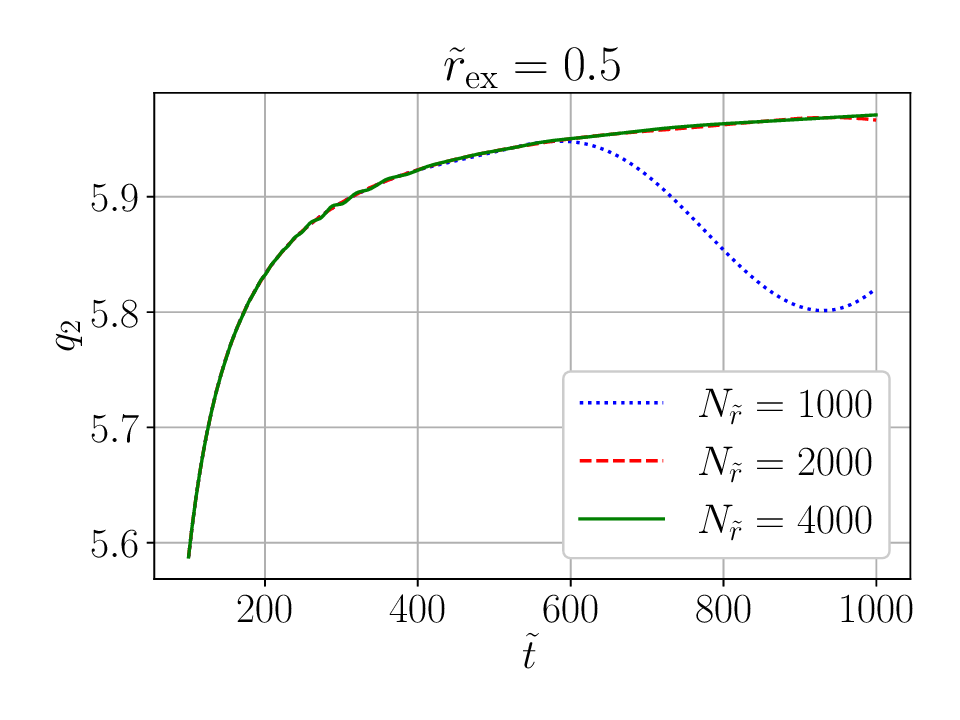}
    \caption{\label{f:lpi_conv_rad_finiter}}
  \end{subfigure}
  \begin{subfigure}[b]{0.5\textwidth}
    \includegraphics[width=\textwidth]{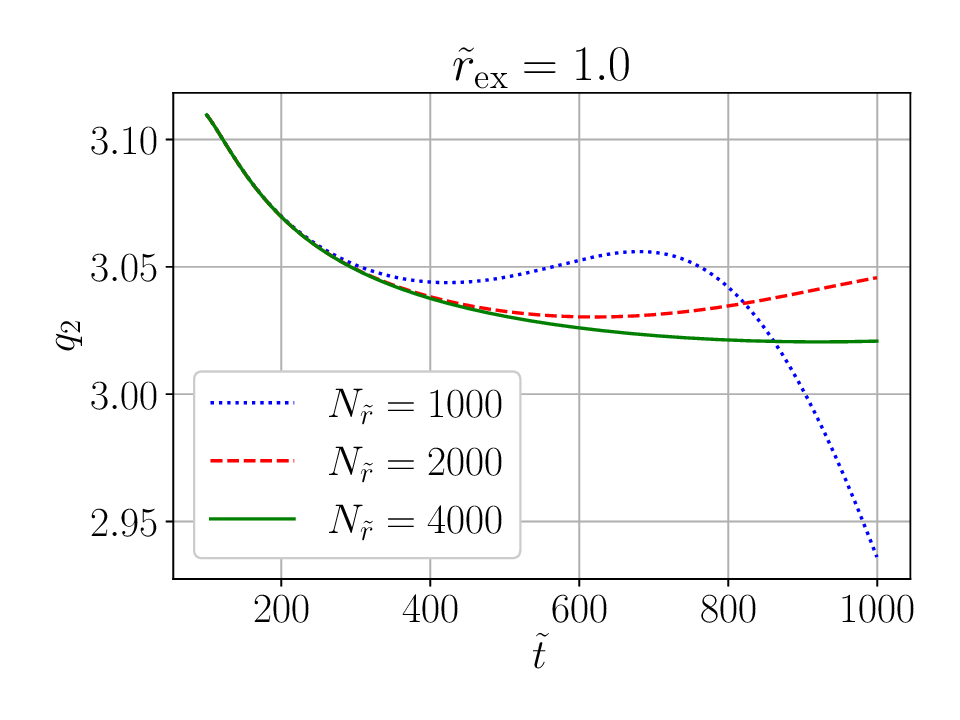}
    \caption{\label{f:lpi_conv_rad_scri}}
  \end{subfigure}
  \begin{subfigure}[b]{0.5\textwidth}
    \includegraphics[width=\textwidth]{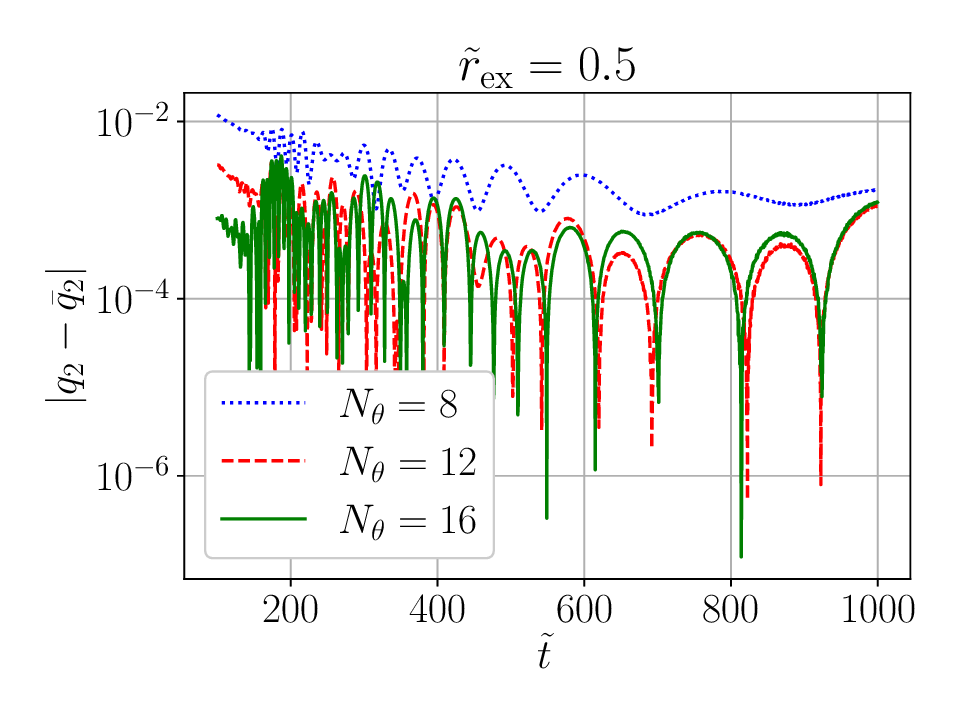}
    \caption{\label{f:lpi_conv_ang_finiter}}
  \end{subfigure}
  \begin{subfigure}[b]{0.5\textwidth}
    \includegraphics[width=\textwidth]{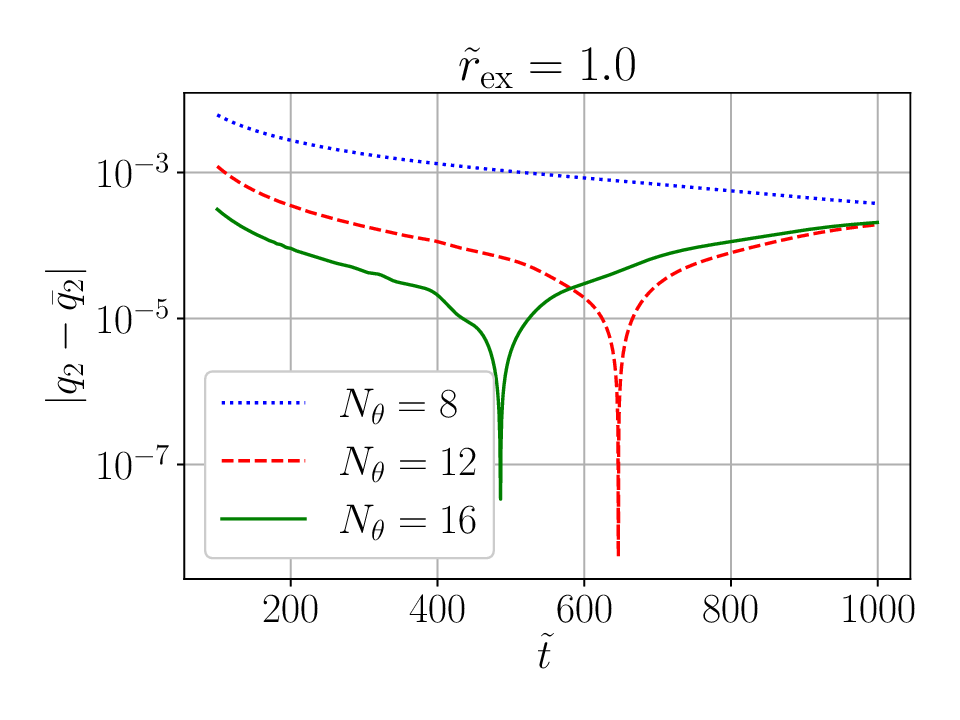}
    \caption{\label{f:lpi_conv_ang_scri}}
  \end{subfigure}
  \caption{\label{f:lpiconv}  
    Convergence test for the local power index in $n=3$ dimensions under axisymmetry with a $p=5$ focusing nonlinearity. 
    Shown is the local power index $q_2$ of the $l=2$ mode extracted at finite radius (a,c) and at \Scri (b,d).
    Static initial data containing two modes with $l=2$ and $l=3$ are chosen, both with $A=12,\, \tr_0 = 0.3$ and $\sigma = 0.07$. 
    (a,b) Varying radial resolutions at fixed angular resolution $N_\theta=12$.
    (c,d) Varying angular resolutions at fixed radial resolution $\Ntr=4000$.
    Here the differences of $q_2$ w.r.t. a reference solution $\bar q_2$ are shown, which was generated using $N_\theta=20$.} 
\end{figure}

In figure \ref{f:lpiconv} we demonstrate convergence of the numerically computed local power index.
From this figure we deduce that increasing the angular resolution beyond $N_\theta=12$ (at fixed radial resolution $\Ntr=4000$) does not lead to more accurate results as far as the local power indices are concerned.
Hence we will always choose $N_\theta=12$ and the highest radial resolution $\Ntr=4000$ in the following.
 
\begin{figure} 
  \centering \includegraphics[width=.5\textwidth]{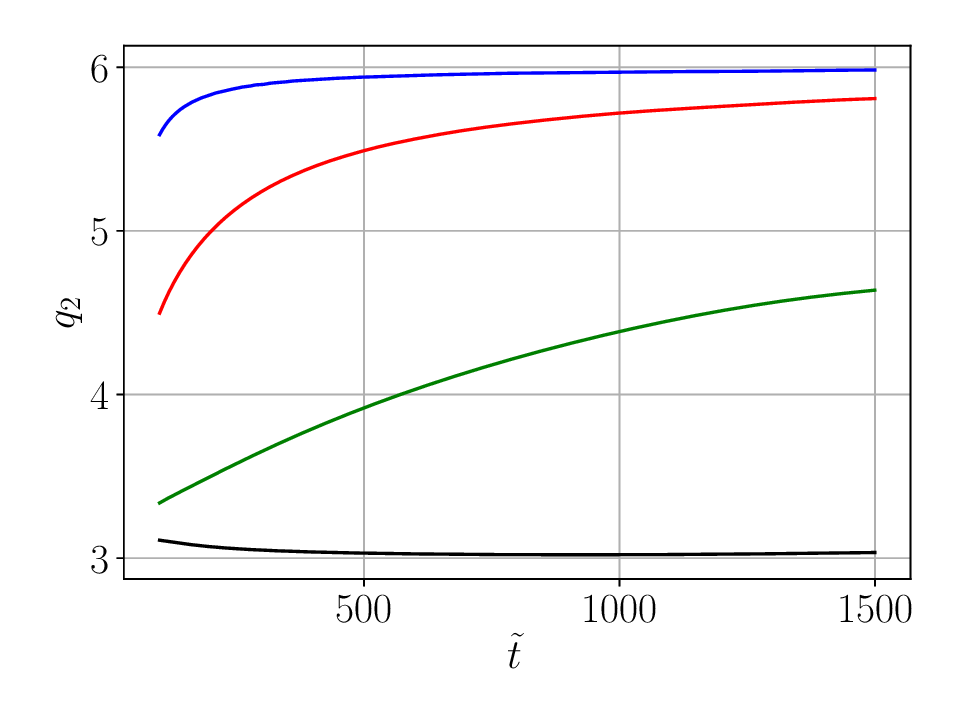}
  \caption{\label{f:lpiradii} 
    Dependence of the local power index on the extraction radius $r_\mathrm{ex}$ in $n=3$ dimensions under axisymmetry with a $p=5$ focusing nonlinearity. 
    Shown is the local power index $q_2$ of the $l=2$ mode extracted at four different radii, from top to bottom: $\tr_\mathrm{ex}=0.5$ (blue), $0.9$ (red), $0.99$ (blue) and $1$ (black).
    Static initial data containing two modes with $l=2$ and $l=3$ are chosen, both with $A=12,\, \tr_0 = 0.3$ and $\sigma = 0.07$. 
    The numerical resolution is $\Ntr=4000, \, N_\theta=12$.
  }
\end{figure}

Next we compare the time evolution of the local power index at different extraction radii in figure \ref{f:lpiradii}. 
We observe that the local power index of a given mode appears to approach the same constant at all finite extraction radii, but in general a different constant at \Scri ($\tr=1$).   

\begin{figure} 
  \begin{subfigure}[b]{0.5\textwidth}
    \includegraphics[width=\textwidth]{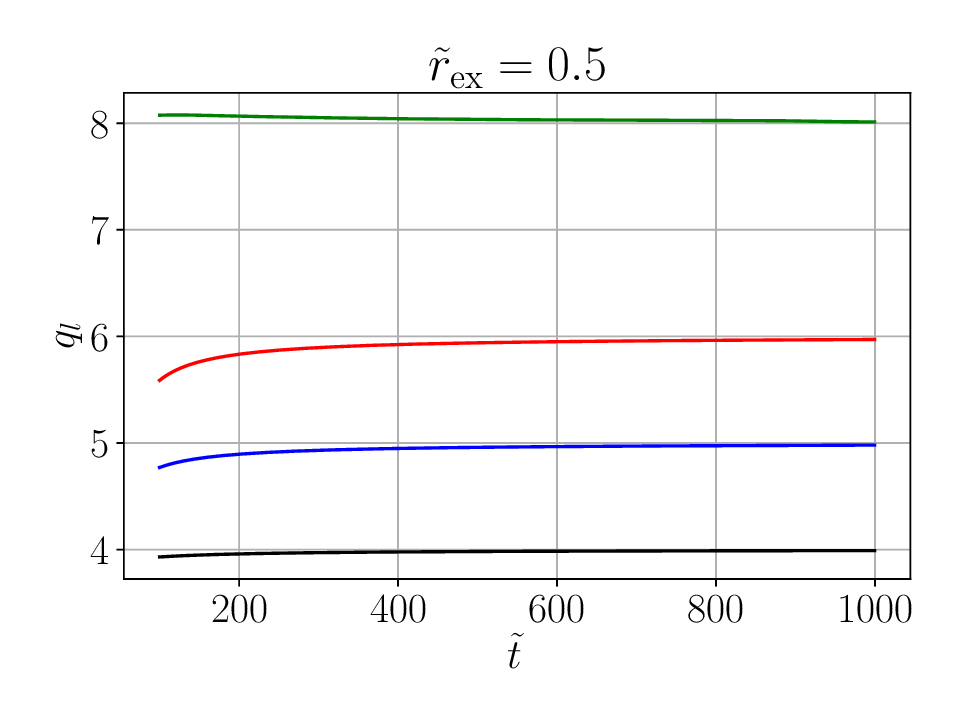}
    \caption{}
  \end{subfigure}
  \begin{subfigure}[b]{0.5\textwidth}
    \includegraphics[width=\textwidth]{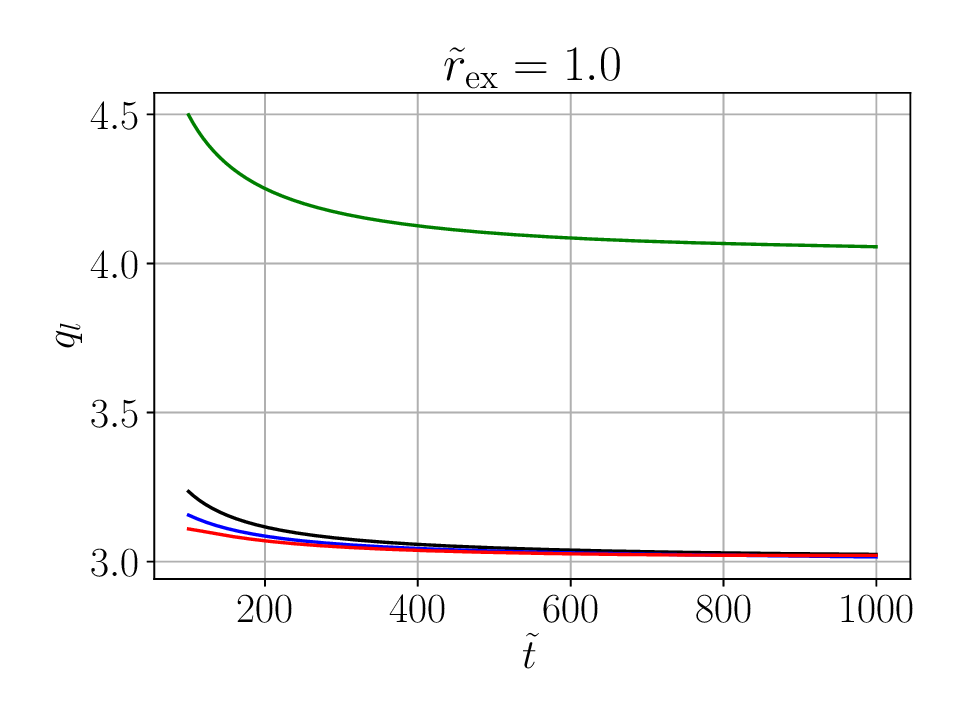}
    \caption{}
  \end{subfigure}
  \caption{\label{f:tailsn3p5} 
    Time evolution of the local power indices $q_l$ in $n=3$ dimensions under axisymmetry with a $p=5$ focusing nonlinearity, extracted (a) at a finite radius and (b) at \Scri.
    From bottom to top in (a): $l=0$ (black), $1$ (blue), $2$ (red) and $3$ (green).
    In (b) the curves for $l=0,1,2$ almost coincide.
    Static initial data containing two modes with $l=2$ and $l=3$ are chosen, both with $A=12,\, \tr_0 = 0.3$ and $\sigma = 0.07$. 
    The numerical resolution is $\Ntr=4000, \, N_\theta=12$.
  }
\end{figure}

\begin{figure} 
  \begin{subfigure}[b]{0.5\textwidth}
    \includegraphics[width=\textwidth]{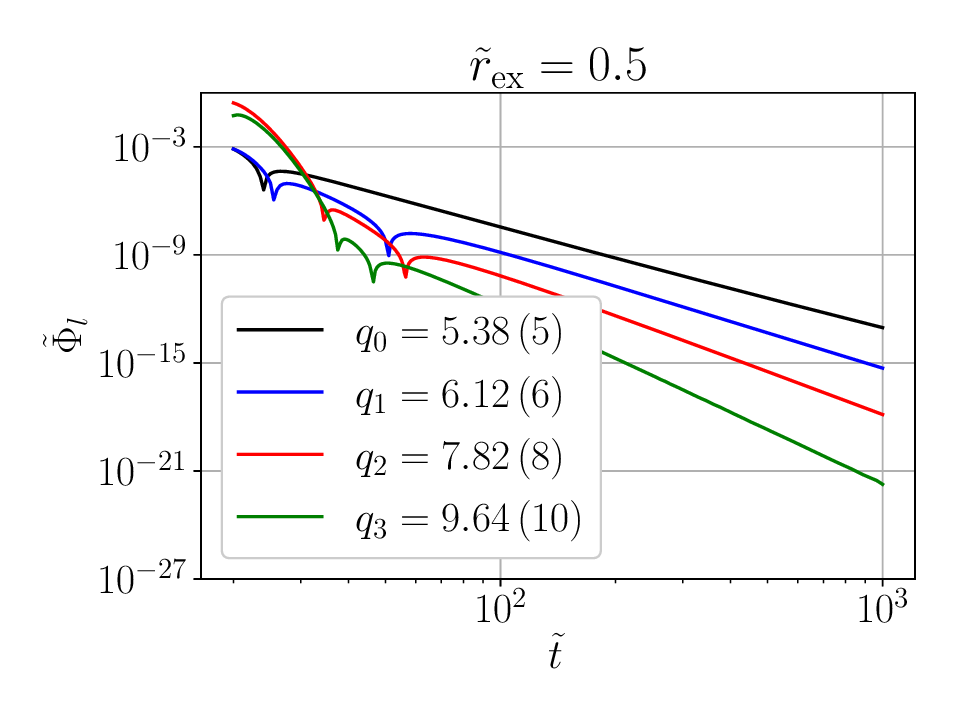}
    \caption{}
  \end{subfigure}
  \begin{subfigure}[b]{0.5\textwidth}
    \includegraphics[width=\textwidth]{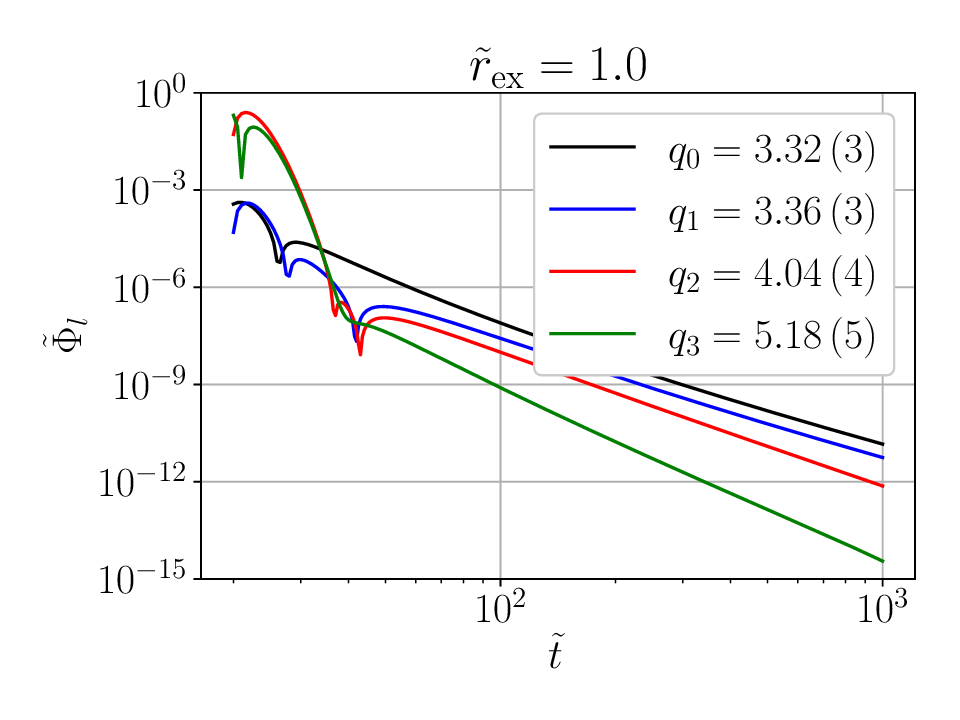}
    \caption{}
  \end{subfigure}
  \caption{\label{f:tailsn5p3} 
    Time evolution of the modes $\tilde\Phi_l$ in $n=5$ dimensions under SO$(4)$ symmetry with a $p=3$ defocusing nonlinearity, extracted (a) at a finite radius and (b) at \Scri.
    The curves refer to $l=0, 1, 2, 3$ from top to bottom.
    The legend shows the fitted decay rates $q_l$ along with their integer estimates (in parentheses).
    Static initial data containing two modes with $l=2$ and $l=3$ are chosen, both with $A=20,\, \tr_0 = 0.3$ and $\sigma = 0.07$. 
    The numerical resolution is $\Ntr=4000, \, N_\theta=12$.
  }
\end{figure}

Finally we investigate the dependence of the decay rates on the spherical harmonic index $l$.
We show results for two selected cases, $n=3,\, p=5$, focusing (figure \ref{f:tailsn3p5}) and $n=5, \, p=3$, defocusing (figure \ref{f:tailsn5p3}).
In the $n=5$ case the higher-$l$ modes are very small.
We had to use higher precision (\texttt{longdouble}) for this run. 
Due to the smallness of the modes, computation of the local power index via \eref{e:lpidef} becomes numerically unstable.
Instead, we have performed a power-law fit to the modes $\tilde \Phi_l$ shown in figure \ref{f:tailsn5p3} in order to obtain the decay rates:
a function of the form $c\cdot \tilde t^{-q_l}$ is fitted by least squares 
(using \texttt{scipy.optimize.curve\_fit}) to each mode $\tilde \Phi_l (\tilde t)$ in the tail phase, i.e. roughly in the interval $\tilde t \in [500,1000]$.

Table \ref{t:decayrates} summarises the asymptotic decay rates $q_l$ found in our numerical evolutions for $n=3$ and $n=5$ and various values of $p$ (subcritical and supercritical).
We have compared focusing and defocusing evolutions and find the same asymptotic decay rates in both cases. 
In addition to the static initial data \eref{e:staticID1}--\eref{e:staticID2} used in the evolutions shown here, we have also tried different initial data taken from the exact linear solutions derived in \ref{s:exact}.
The observed decay rates are identical. 
The static data have the advantage that they have an initially outgoing component so the tail forms earlier and with a higher amplitude than for the data from the linear solutions, which are initially ingoing.

\begin{table}
  \[
  \begin{array}{c||c|c|c|c|c}
    n=3 & p=3 & p=4 & p=5 & p=6 & p=7 \\\hhline{=||=|=|=|=|=}
    l=0 & 2|1 & 3|2 & 4|3 & 5|4 & 6|5 \\\hline
    l=1 & 4|2 & 4|2 & 5|3 & 6|4 & 7|5 \\\hline
    l=2 & 6|3 & 6|3 & 6|3 & 7|4 & 8|5 \\\hline
    l=3 & 8|4 & 8|4 & 8|4 & 8?|4 & 9?|5? 
  \end{array}
  \]
  \[ 
  \begin{array}{c||c|c}
    n=5 & p=2 & p=3 \\\hhline{=||=|=}
    l=0 & 4|2 & 5?|3? \\\hline
    l=1 & 6|3 & 6|3? \\\hline
    l=2 & 8|4 & 8|4 \\\hline
    l=3 & 10|5 & 10?|5
  \end{array}
  \]
  \caption{\label{t:decayrates}
    Asymptotic decay rates $q_l$ for $n=3$ and $n=5$ under SO$(n-1)$ symmetry, for various values of the power $p$ of the nonlinearity, determined from numerical evolutions.
    For each $l$, the values of $q_l$ at a finite extraction radius (left) and at \Scri (right) are shown. 
    Values marked ? are uncertain.
  }
\end{table}


\section{Conclusion and discussion}
\label{s:conclusion}

This paper is concerned with the nonlinear wave equation 
\begin{equation} \label{e:nlw3}
  \Box \Phi := -\partial_t^2 \Phi + \Delta \Phi =  \mu |\Phi|^{p-1} \Phi, \quad \Phi : \R \times \R^n \to \R
\end{equation}
with $p>1$, $\mu=\pm1$ in $n\geqslant 3$ spatial dimensions.
Unlike previous studies, we do not assume radial symmetry. 
In $n=3$ we do not assume any symmetries, and in higher dimensions we impose an SO$(n-1)$ symmetry so that there is one effective angular coordinate. 

We introduce a foliation of Minkowski spacetime into hyperboloidal slices of constant mean curvature, combined with a conformal compactification. 
This avoids the need for artificial timelike boundaries and allows us to numerically construct the solution in the entire future of the initial hyperboloidal slice, all the way up to future null infinity \Scri, where future-directed null gedodesics end.

Our numerical approach combines a fourth-order finite difference method in the radial coordinate with a Fourier pseudo-spectral method in the angular coordinate(s).
We construct exact solutions to the linear wave equation and demonstrate fourth-order convergence of our code as the radial resolution is increased.

Unlike on standard slices approaching spatial infinity, the energy of the scalar field on hyperboloidal slices is not conserved.
It decreases due to the energy flux at \Scri.
We derive this energy balance and show that it is well satisfied in our numerical evolutions.

The main result of this paper concerns the late-time power-law decay of the solutions.
We expand the evolved field into spherical harmonics and determine decay exponents of the individual modes for various values of the exponent $p$ of the nonlinearity (subcritical, critical and supercritical).
The decay exponents appear to be the same for focusing and defocusing nonlinearities ($\mu=\pm 1$).
They are independent of the extraction radius as long as this radius is finite, but at \Scri different (smaller) decay rates are found.
In three dimensions the decay exponents are found to be independent of the azimuthal spherical harmonic index $m$, which justifies our assumption of axisymmetry or more generally SO$(n-1)$ symmetry.
Our results in dimensions $n=3$ and $n=5$ are summarised in table \ref{t:decayrates}.

At a finite radius, we would expect the scalar field to show the same decay in standard coordinates $t, r$ as in hyperboloidal coordinates $\tilde t, \tr$.
This leads us to the following
\begin{conj}
  Consider the nonlinear wave equation \eref{e:nlw3} in $n=3$ spatial dimensions with $p\in\mathbb{N}$, $p\geqslant 3$.
  Expanding the field in spherical harmonics
  \begin{equation}
    \Phi(t, r, \theta, \varphi) = \sum_{l=0}^\infty \sum_{m=-l}^l \Phi_{lm} (t, r) \, Y_{l m}(\theta, \varphi),
  \end{equation}
  the modes decay asymptotically (as $t\to\infty$) as 
  \begin{equation}
    \Phi_{lm}(t, r) \sim t^{-q_l}, \qquad  q_l = \max(l+p-1, \, 2l+2)
  \end{equation}
  at any fixed finite radius $r$.
  On slices $\Sigma_{\tilde t}$ approaching future null infinity \Scri, where $t, r \to\infty$, $t - r \to \tilde t$ and the conformal factor $\Omega(r) = \Or(r^{-1})$, the conformally rescaled scalar field $\tilde \Phi =\Omega(r)^{(1-n)/2} \, \Phi$ decays at \Scri asymptotically as 
  \begin{equation}
    \tilde \Phi_{lm}(\tilde t) \sim \tilde t^{-\tilde q_l}, \qquad \tilde q_l = \max(p-2, \, l+1).
  \end{equation}
\end{conj}

Indeed, this is consistent with \cite{Szpak2009}, where the authors proved using perturbative methods that in $n=3$ spatial dimensions under spherical symmetry ($l=0$) and $p\geqslant 3$ the solution decays asymptotically as $t^{-p+1}$ at a finite radius.
It would be very interesting to try and prove the above conjecture for $l>0$.

In dimension $n=5$, a similar form of the decay exponents consistent with table \ref{t:decayrates} would be 
\begin{equation}
  q_l = \max(l+p+2, \, 2l+4), \qquad \tilde q_l = \max(p, \, l+2)
\end{equation}
for the decay at a finite radius and at \Scri, respectively.
More data would be needed to corroborate such a conjecture for $n=5$ though. 
This proves difficult numerically because the field decays very rapidly for higher exponents $p$ of the nonlinearity.
Higher than \texttt{longdouble} precision would be needed.

The observation that the solution decays more slowly at \Scri than at any finite radius means that the solution will develop an increasingly steep radial gradient near $\tr=1$ as time progresses. 
For the simulations shown in this paper and for the runtimes needed to safely determine the decay exponents ($t_{\max} \approx 1000$), we did not observe a significant loss of accuracy or failure of convergence, cf. figure \ref{f:lpiconv}.
For longer runtimes, some form of grid adaptivity (non-uniform grids or adaptive mesh refinement) will eventually be needed in order to resolve this feature.

Apart from tails, there are various other aspects of nonlinear wave equations that could be investigated using our hyperboloidal evolution code, in particular the nature of singularity formation (blow-up) and the threshold between blow-up and scattering.
For example, this threshold was investigated numerically for the focusing cubic ($p=3$) wave equation in the subcritical dimension $n=3$ in radial symmetry in \cite{Bizon2009} (see also a similar study for the nonlinear Klein-Gordon equation \cite{Donninger2011}).
In \cite{Glogic2020} the focusing cubic wave equation was considered in the supercritical dimensions $n=5$ and $n=7$ and again threshold solutions were identified.
It will be interesting to analyse these situations beyond radial symmetry.


\appendix
 
\section{Exact solutions, initial data}
\label{s:exact}

Here we work out exact solutions to the linear wave equation, which are used to test the convergence of the code in section \ref{s:convtest}.
We first assume dimension $n=3$ without symmetries in \ref{s:exactn3}; higher dimensions will be treated in \ref{s:exactn5}.
The solutions are constructed in standard spherical polar coordinates first;
finally in \ref{s:exacttrafo} they are transformed to hyperboloidal coordinates.


\subsection{$n=3$ without symmetries} \label{s:exactn3}

The solution to the linear wave equation 
\begin{equation} 
  \Box \Phi = 0
\end{equation}
in standard spherical polar coordinates can be decomposed into spherical harmonics,
\begin{equation} \label{e:spherharmexpansion}
  \Phi = \sum_{l=0}^\infty \sum_{m=-l}^l\Phi_{lm} (t, r) \, Y_{l m}(\theta, \varphi).
\end{equation}
The spherical harmonics are eigenfunctions of the Laplace-Beltrami operator on $S^2$,
\begin{equation}
  \mathring{\Delta}^{(2)} \, Y_{l m}(\theta, \varphi) = -l(l+1) Y_{l m}(\theta, \varphi).
\end{equation}
They have the form
\begin{equation}
  Y_{l m}(\theta,\varphi) = N_{lm} P_{l m} (\theta) \, \rme^{\rmi m \varphi},
\end{equation}
where $P_{l m} (\theta)$ are the associated Legendre functions and $N_{lm}$ are normalisation constants chosen such that 
\begin{equation}
  \fl \int_{S^{2}} Y_{lm} Y_{l'm'} dS^{(2)} = 
  \int_0^{2\pi} \int_0^\pi Y_{lm}(\theta,\varphi) Y_{l'm'}^*(\theta,\varphi) \sin\theta \,\rmd\theta \, \rmd\varphi = \delta_{ll'} \delta_{mm'}.
\end{equation}
The (complex) spherical harmonics are implemented in 
\texttt{scipy.special.sph\_harm\_y}.
A suitable real basis $\{ \hat Y_{lm} \}$ of spherical harmonics is obtained by taking 
\begin{equation} \label{e:real_sph_harm}
  \hat Y_{lm} := \left\{
  \begin{array}{ll}
    \sqrt{2} (-1)^m \Im Y_{l|m|}, & m < 0, \smallskip\\
    Y_{l0}, & m=0, \smallskip\\ 
    \sqrt{2} (-1)^m \Re Y_{lm}, & m > 0.
  \end{array}
  \right.
\end{equation}

The wave equation implies that the radial functions $\Phi_{lm}(t,r)$ in \eref {e:spherharmexpansion} obey an equation of Euler-Poisson-Darboux type,
\begin{equation}  \label{e:epd3}
  -\Phi_{lm,tt} + \Phi_{lm,rr} + \frac{2}{r} \Phi_{lm,r} - \frac{l(l+1)}{r^2} \Phi_{lm} = 0.
\end{equation}
Since this equation does not depend on $m$ explicitly, we leave out the index $m$ on $\Phi_{lm}$ in the following.
Solutions have the form
\begin{equation}
  \Phi_l^\pm(t, r) = \sum_{k=0}^l c_k \, r^{-k-1} F^{(l-k)}(r \pm t),
\end{equation}
where $F$ is an arbitrary mode function, and $-$ refers to an outgoing and $+$ to an ingoing solution.
The coefficients $c_k$ can be determined by recursion.
We provide explicit solutions for $l=0,1,2$:
\begin{eqnarray}
  \label{e:Phi_rad0_n3}
  \Phi_0^\pm(t,r) &=& r^{-1} F(r\pm t), \\
  \Phi_1^\pm(t,r) &=& r^{-1} F'(r\pm t) - r^{-2} F(r\pm t),\\
  \Phi_2^\pm(t,r) &=& r^{-1} F''(r\pm t) - 3 r^{-2} F'(r\pm t) + 3 r^{-3} F(r\pm t).
  \label{e:Phi_rad2_n3}
\end{eqnarray}
We observe that $\Phi_l^- = \Or(r^{-1})$ as \Scri is approached ($r\to\infty$, $r-t$ finite).

If we take the mode function $F$ to be odd and form the linear combination $\Phi_l^+ + \Phi_l^-$, the result turns out to be regular at the origin $r = 0$.

 
\subsection{$n\geqslant 3$ with $\mathrm{SO}(n-1)$ symmetry}
\label{s:exactn5}

In this case we expand
\begin{equation}
  \Phi = \sum_{l=0}^\infty \Phi_{l} (t, r) \, Y_l(\theta),
\end{equation}
where the real spherical harmonics $Y_l(\theta)$ are eigenfunctions of the Laplace-Beltrami operator on $S^{(n-1)}$,
\begin{equation}
  \mathring{\Delta}^{(n-1)} \, Y_l(\theta) = l(2-n-l) Y_l(\theta).
\end{equation}
We normalise them such that
\begin{equation}
  \int_{S^{n-1}} Y_l Y_{l'} dS^{(n-1)} = 
  A^{(n-2)} \int_0^\pi Y_l(\theta) Y_{l'}(\theta) \sin^{n-2}\theta \,\rmd\theta = \delta_{ll'},
\end{equation}
where the area $A^{(n-2)}$ of $S^{n-2}$ is given in \eref{e:areaS}.
The mode equation now reads
\begin{equation}
  -\Phi_{l,tt} + \Phi_{l,rr} + \frac{n-1}{r} \Phi_{l,r} + \frac{l(2-n-l)}{r^2} \Phi_l = 0.
\end{equation}
Note that for $n=3$ we recover \eref{e:epd3}.

In the following we state explicit solutions for $n=5$.
The first three spherical harmonics are
\begin{equation} \fl
  Y_0(\theta) = \frac{\sqrt{3}}{2\sqrt{2}\pi}, \quad
  Y_1(\theta) = \frac{\sqrt{15}}{2\sqrt{2}\pi} \cos\theta, \quad 
  Y_2(\theta) = \frac{\sqrt{21}}{8\pi} (5\cos^2\theta - 1)
\end{equation}
and the corresponding radial solutions are 
\begin{eqnarray}
  \label{e:Phi_rad0_n5}
  \fl \Phi_0^\pm( t, r) &=& r^{-2} F'(r\pm  t) - r^{-3} F(r\pm  t),  \\
  \fl \Phi_1^\pm( t, r) &=& r^{-2} F''(r\pm  t) - 3 r^{-3} F'(r\pm  t) + 3 r^{-4} F(r\pm  t), \\
  \label{e:Phi_rad2_n5}
  \fl \Phi_2^\pm( t, r) &=& r^{-2} F'''(r\pm  t) - 6 r^{-3} F''(r\pm  t) + 15 r^{-4} F'(r\pm  t) -15 r^{-5} F(r\pm  t).
\end{eqnarray}
We observe that now $\Phi_l^- = \Or(r^{-2})$ as \Scri is approached.


\subsection{Transformation to hyperboloidal coordinates}
\label{s:exacttrafo}

From the conformal coordinates $\tilde t, \tr$ used in the code, we first compute the physical coordinates  
\begin{equation}
  r = \frac{2n\tr}{C(1-\tr^2)}, \quad t = \tilde t + \sqrt{\frac{n^2}{C^2} + r^2}
\end{equation}
so we can compute $\Phi(t,r)$ as constructed above.
We then form 
\begin{equation} \label{e:PhifromPhi_}
  \tilde \Phi = \Omega^{(1-n)/2} \Phi,
\end{equation}
where the conformal factor is given by \ref{e:Omega}.

The multiplication by a negative power of $\Omega$ in \eref{e:PhifromPhi_} may seem problematic at \Scri ($r=1$), where $\Omega=0$. 
However, from the explicit solutions given in \ref{s:exactn3} and \ref{s:exactn5}, we see that $\Phi = \Or(r^{(1-n)/2})$,
and $\Omega = \tr/r = \Or(r^{-1})$, so 
$\tilde \Phi$ in \eref{e:PhifromPhi_} is in fact regular at \Scri.

For the initial data we also need to compute the field $\tilde \Pi$.
By its definition, we have 
\begin{equation}
  \fl \tilde \Pi = \tilde \alpha^{-1} (\tilde \Phi_{,\tilde t} - \tilde \beta^{\tr}\Phi_{,\tr})
    = \tilde \alpha^{-1} \Omega^{(1-n)/2} \left[\Phi_{,\tilde t} - \tilde \beta^{\tr} \left(\Phi_{,\tr} + \frac{1-n}{2} (\ln\Omega)_{,\tr} \Phi \right)\right]
\end{equation}
with $\Omega$, $\tilde \beta^r$ and $\tilde \alpha$ given in \eref{e:Omega}--\eref{e:lapse}.
Finally we need to express the derivatives of $\Phi$ in terms of derivatives w.r.t. $r$ and $t$:
\begin{equation}
  \Phi_{,\tr} = \frac{\rmd r}{\rmd \tr} \Phi_{,r} + \frac{\partial t}{\partial \tr} \Phi_{,t}, \qquad \Phi_{,\tilde t} = \Phi_{,t},
\end{equation}
where
\begin{equation}
   \frac{\rmd r}{\rmd \tr} = \frac{2n(1+\tr^2)}{C(1-\tr^2)^2}, \qquad
   \frac{\partial t}{\partial \tr} = \frac{4n\tr}{C(1-\tr^2)^2}.
\end{equation}
The quantities $\Phi_{,r}$ and $\Phi_{,t}$ can be computed directly from the radial solutions \eref{e:Phi_rad0_n3}--\eref{e:Phi_rad2_n3} or \eref{e:Phi_rad0_n5}--\eref{e:Phi_rad2_n5}.


\section*{References}


\providecommand{\newblock}{}

\end{document}